\documentclass{article}

\usepackage[preprint]{neurips_2023}
\usepackage[utf8]{inputenc} % allow utf-8 input
\usepackage[T1]{fontenc}    % use 8-bit T1 fonts

\usepackage[hyphens]{url}           % added
\usepackage{hyperref}               % added
\usepackage[hyphenbreaks]{breakurl} % added

\usepackage{booktabs}       % professional-quality tables
\usepackage{amsfonts}       % blackboard math symbols
\usepackage{nicefrac}       % compact symbols for 1/2, etc.
\usepackage{microtype}      % microtypography
\usepackage{xcolor}         % colors
\usepackage[normalem]{ulem}

\usepackage{amsmath, amsfonts, amsthm, amssymb}
\usepackage{color, bm, graphicx} 
\usepackage{booktabs, multicol, multirow}
\usepackage{xspace}

\usepackage[utf8]{inputenc} 
\usepackage[english]{babel} 
\usepackage{natbib} 
\bibpunct[, ]{(}{)}{,}{a}{}{,}

\usepackage{verbatim} 
\usepackage{array} 
\usepackage{lscape} 
\usepackage{adjustbox} 
\usepackage{booktabs}
\usepackage{multirow}
\usepackage{makecell} 
\usepackage{caption}
\usepackage{threeparttable}
\usepackage{xcolor}
\usepackage{algorithm}
\usepackage{algpseudocode}

\usepackage{endnotes}
\usepackage{float}

\newtheorem{lemma}{Lemma}
\newtheorem{proposition}{Proposition}

\newtheorem{assumption}{Assumption}

\def\Expect{{\mathbb E}}
\def\Prob{{\mathbb P}}

\def\subto{{\rm subject \mbox{   }\rm to}}

\newcommand{\remove}[1]{}

\newcommand{\V}[1]{{{\boldsymbol #1}}}

\newcommand{\x}{\V{x}}

\newcommand{\Vzeta}{\V{\zeta}}

\newcommand{\X}{\mathcal{X}}

\newcommand{\1}{\V{1}}

\newcommand{\Dist}{\mathcal{D}}

\newcommand{\mymbox}[1]{\mbox{\scriptsize #1}}
\renewcommand{\Re}{\mathbb{R}}
\newcommand{\quoteIt}[1]{``#1''}

\newcommand{\Vxi}{{\V{\xi}}}
\newcommand{\hatFC}{\tilde{F}}

\newcommand{\comments}[1]{}
\newcommand{\removed}[1]{}
\newcommand{\EDmodified}[1]{{#1}}
\newcommand{\AGmodified}[1]{{#1}}

\newcommand{\EDcomments}[1]{{}}
\newcommand{\AGcomments}[1]{{}}
\newcommand{\MPcomments}[1]{{}}

\newcommand{\coeffPresc}{{\mathcal{P}}}
\newcommand{\xSAA}{{\hat{\x}}}
\newcommand{\ratioPresc}{{\mathcal{V}}}
\newcommand{\coeffDet}{{R^2}}
\newcommand{\policySet}{{\mathcal{H}}}
\newcommand{\DistCVaR}{{\bar{\Dist}}}

\def\min{\mathop{\rm min}}
\def\max{\mathop{\rm max}}
\def\argmax{\mathop{\rm argmax}}
\def\argmin{\mathop{\rm argmin}}
\newcommand{\CVaR}{\mbox{CVaR}}
\newcommand{\pushToAppendix}[1]{{}}

\usepackage{subcaption}
\usepackage{appendix}
\title{Robust Data-driven Prescriptiveness Optimization}

\author{%
  Mehran Poursoltani\\
  Desautels Faculty of Management\\
  McGill University\\
  Montréal, Québec, Canada \\
  \texttt{mehran.poursoltani@mcgill.ca} \\
  \And
  Erick Delage \\
  GERAD \& Department of Decision Sciences \\
  HEC Montréal \\
  Montréal, Québec, Canada \\
  \texttt{erick.delage@hec.ca} \\
  \AND
  Angelos Georghiou \\
  Department of Business and Public Administration \\
  University of Cyprus \\
  Nicosia, Cyprus \\
  \texttt{georghiou.angelos@ucy.ac.cy}
}
\begin{document}
\maketitle
\begin{abstract}
    The abundance of data has led to the emergence of a variety of optimization techniques that attempt to leverage available side information to provide more anticipative decisions. The wide range of methods and contexts of application have motivated the design of a universal unitless measure of performance known as the coefficient of prescriptiveness. This coefficient was designed to quantify both the quality of contextual decisions compared to a reference one and the prescriptive power of side information. To identify policies that maximize the former in a data-driven context, this paper introduces a distributionally robust \EDmodified{contextual} optimization model where the coefficient of prescriptiveness substitutes for the classical empirical risk minimization objective. We present a bisection algorithm to solve this model, which relies on solving a series of linear programs when the distributional ambiguity set has an appropriate nested form and polyhedral structure. Studying a contextual shortest path problem, we evaluate the robustness of the resulting policies against alternative methods when the out-of-sample dataset is subject to varying amounts of distribution shift.
\end{abstract}

\section{Introduction}\label{sec:intro}

Stochastic programming is perceived as one of the fundamental methods devised for decision-making under uncertainty (see \cite{shapiro2021lectures} and \cite{birge2011introduction}). Given a cost function $h(\x,\V{\xi})$ that depends on a decision $\x\in\Re^{n_\x}$ and a random vector $\V{\xi}\in\Re^{n_\V{\xi}}$, the stochastic programming (SP) problem is defined as
\begin{equation}
(SP)\;\;\;\x^*\in\arg\min_{\x\in\X}\;\;\Expect_{F}[h(\x,\V{\xi})],\label{eq:SP}
\end{equation}
where $\X$ is a convex feasible set, $h(\x,\V{\xi})$ is a cost function that is assumed convex in $\x$ for all $\Vxi$, and $\V{\xi}$ is assumed to be drawn from the distribution $F$. The solution methods for this problem mainly rely on either assuming a priori distribution for $F$ or exploiting a set of independent and identically distributed observations. In the latter case, a set of i.i.d observations of the random vector $\V{\xi}$ denoted by $\mathcal{S}:=\{\V{\xi_i}\}_{i=1}^N$ can be used to formulate  the following sample average approximation problem:
\begin{equation}
(SAA)\;\;\;\x^*\in\arg\min_{\x\in\X}\;\;\frac{1}{N}\sum_{i=1}^N h(\x,\V{\xi_i}),\label{eq:SP_data}
\end{equation}
where we assume a uniform distribution over the observed data. Recently, the availability of large datasets has played a critical role in redirecting the optimization methods devised for decision-making under uncertainty towards taking advantage of so-called \quoteIt{side information} or \quoteIt{covariates}. This paradigm encourages decision-makers to benefit from the available data beyond the desired random variables to make more anticipative decisions. For instance, a portfolio manager who optimizes her investments in the stock market may consider a variety of available micro and macroeconomic indicators as side information to make more anticipative decisions (see \cite{10.1093/rfs/hhp003} and \cite{matte2020}), \AGmodified{while a traffic path planner can utilize side information like time of day, weather status and holiday/work day to find the best route through the city (see \cite{BertsimasKallus})}. This gives rise to the following {\color{black}contextual} stochastic optimization (CSO) problem:
\begin{equation}
(CSO)\;\;\;\x^*(\Vzeta)\in\arg\min_{\x\in\X}\;\;\Expect_{F}[h(\x,\V{\xi})|\Vzeta],\label{eq:CSO}
\end{equation}
where $\Vzeta\in\Re^{n_{\Vzeta}}$ denotes the given vector of \quoteIt{covariates}, or so-called \quoteIt{features}. In this case, any observed random vector $\V{\xi_i}$ is accompanied by a vector of covariates $\Vzeta_i\in\Re^{n_\Vzeta}$. The difficulty of this problem shows up when the conditional probability distribution function $F_{\V{\xi}|\Vzeta}$ is unknown, and only a set of i.i.d observations $\mathcal{T}:=\{(\Vzeta_i,\V{\xi_i})\}_{i=1}^N$ is available. In this case, a data-driven variant of the CSO problem can be written as
\begin{equation}
\x^*(\Vzeta)\in\arg\min_{\x\in\X} \Expect_{\hat{F}_{\xi|\zeta}}[h(\x,\V{\xi})],\label{eq:CSO_data}
\end{equation}
where $\hat{F}_{\xi|\zeta}$ is a conditional probability model for $\Vxi$ given $\Vzeta$ inferred from the available data, e.g. by training a random forest \citep{breiman2001random}, or estimated via kernel density estimation \citep{ban2019big}. \AGmodified{To deal with  possible overfitting in the presence of limited data or possible distribution shifts due to unexpected events, one can formulate a distributionally robust {\color{black}contextual} stochastic optimization (DRCSO) model, which in general, takes the following form
\begin{equation}
(DRCSO)\,\,\,\x^*(\Vzeta)\in\arg\min_{\x\in\X}\sup_{F\in\mathcal{D}}\;\;\Expect_{F}[h(\x,\V{\xi})|\Vzeta]\label{eq:DRCSO}
\end{equation}
where $\mathcal{D}$ is the ambiguity set containing the set of admissible distributions (see \cite{DuchiDRLossesLatentCovMix,bertsimas2022bootstrap,kannan2020residuals,delageDRCSO,estebanperez2022,hanasusantoDRCSO}, and literature within). 
}

Recently, \cite{BertsimasKallus} proposed to compare the performance of different CSO (or DRCSO) approaches, by measuring the \quoteIt{coefficient of prescriptiveness}, defined as:
\begin{equation}
    \coeffPresc_F(\x(\cdot)):= 1- \frac{\Expect_F[h(\x(\Vzeta),\V{\xi})]-\Expect_F[\min_{\x'\in\X}h(\x',\V{\xi})]}{\Expect_F[h(\xSAA,\V{\xi})]-\Expect_F[\min_{\x'\in\X}h(\x',\V{\xi})]},
\end{equation}
where $\xSAA:= \argmin_{\x}\Expect_{\hat{F}}[h(\x,\V{\xi})]$ with $\hat{F}$ as the empirical distribution that puts equal weights on each observed data point $\{\V{\xi_i}\}_{i=1}^N$ (i.e. the solution of SAA). The idea behind the coefficient of prescriptiveness is that it measures the performance of a given policy $\bm x(\bm \zeta)$ relative to the constant decision $\xSAA$ which is agnostic to the side information $\bm \zeta$, and to the fully anticipative policy which achieves the progressive optimal value of $\Expect_F[\min_{\x'\in\X}h(\x',\V{\xi})]$. It is easy to see that  a high value of $\mathcal{P}_F$ indicates that the policy can leverage the contextual information of $\bm\zeta$ with $\mathcal{P}_F = 1$ indicating that the policy is achieving the fully anticipative performance in terms of $\bm\xi$. In contrast, a low value of $\mathcal{P}_F$ indicates that the policy is not able to exploit (or even is misled by) the available information. \AGmodified{This behavior is reminiscent of $R^2$, the ``coefficient of determination", typically used in the context of predictive models, a connection which we discuss in the next section.}

Following the introduction of the coefficient of prescriptiveness, this metric has been employed in several pieces of research to demonstrate the potential of proposed data-driven policies for leveraging the available side information. One can refer to \cite{bertsimas2016inventory} for such a comparison in the context of inventory management, \cite{stratigakos2022prescriptive} for energy trading, \cite{Notz2022} for flexible capacity planning, and \cite{StochasticOptimizationForests} for shortest path and portfolio optimization problems. We note that, in the current literature, $\mathcal{P}_F$ is only used as a benchmark metric for assessing the performance of policies computed using different approaches, e.g., in \cite{BertsimasKallus}, the metric compares policies computed (amongst others) using CSO where the conditional probability is estimated by random forests and kernel density estimation. Given its prevalence as a performance measure, it is natural to question whether it is possible and useful to directly optimize the coefficient of prescriptiveness.

While one can show that maximizing $\coeffPresc_F$ reduces to solving the CSO problem, one may wonder how the $\coeffPresc_F$ measure should be robustified in order to improve out-of-sample performance. In this work, we introduce for the first time a distributionally robust version of $\coeffPresc_F$. We establish connections to other models in the literature and present an efficient algorithm to maximize it when the conditional probability model is discrete (such as with a random forest or with a Kernel density estimator).
The rest of the paper is organized as follows. Section \ref{sec:motivationCoeffPresc} motivates the optimization of the coefficient of prescriptiveness by explicating its relationship to the coefficient of determination in the field of statistics. Section \ref{sec:Model} introduces a robust data-driven prescriptiveness optimization model that can be used to maximize a distributionally robust version of the coefficient of prescriptiveness. We reformulate this problem as a convex optimization problem that can reduce to a linear program when the ambiguity set takes the form of a so-called \quoteIt{nested Conditional Value-at-Risk (CVaR) set}. A bisection method is proposed to solve the latter, as well as an acceleration scheme; finally, Section \ref{sec:Experiments} presents the numerical experiments, where we evaluate the robustness of the resulting policies against benchmark ones in a shortest path problem when the out-of-sample dataset confronts a distribution shift. All proofs are relegated in Appendix~\ref{appendix:proofs}.

\section{Motivation for optimizing $\coeffPresc$ and its robustification}\label{sec:motivationCoeffPresc}

As argued in \cite{BertsimasKallus}, in the context of predictive models, where one wishes to     predict the value of $\xi\in\Re$ based on a list of covariates $\Vzeta$ using a statistical model $f:\Re^{n_\zeta}\rightarrow \Re$, one popular metric that is employed takes the form of the so-called \quoteIt{coefficient of determination}:
\[\coeffDet(f(\cdot)):= 1- \frac{\Expect_{\hat{F}}[(f(
\Vzeta)-\xi)^2]}{\Expect_{\hat{F}}[(\hat{\xi}-\xi)^2]},\]
where $\hat{\xi}:=\Expect_{\hat{F}}[\xi]$ is the empirical mean of $\xi$ in the data set {\color{black}and $\hat F$ is the empirical joint distribution of $(\bm\zeta,\xi)$}. The popularity of $\coeffDet{}$ compared to mean squared error as a measure of performance can be partially attributed to being unitless. It is upper bounded by 1, with a value closer to $1$, indicating that most of the variation of $\xi$ can be modeled using $f(\cdot)$. On the flip side, when strictly smaller than 0, its absolute value measures the percentage of additional variations that are introduced by the predictive model, thus indicating a degradation of predictive power when compared to the simple sample average $\hat{\xi}$.

The coefficient of prescriptiveness can be viewed as an attempt to introduce an analogous measure in the {\color{black}contextual} optimization setting. More specifically, it reduces to $\coeffDet$ when $n_{x}=1$ and $h(x,\xi):=(x-\xi)^2$, namely:
\[\coeffPresc_{\hat{F}}(x(\cdot))=1-\frac{\Expect_{\hat{F}}[(x(\Vzeta)-\xi)^2)]-\Expect_{\hat{F}}[\min_{x'}(x'-\xi)^2]}{\Expect_{\hat{F}}[(\hat{x}-\xi)^2]-\Expect_{\hat{F}}[\min_{x'}(x'-\xi)^2]}=\coeffDet(x(\cdot)),\]
% \begin{equation*}
%     \begin{aligned}
%         \coeffPresc_{\hat{F}}(x(\cdot))&=1-\frac{\Expect_{\hat{F}}[(x(\Vzeta)-\xi)^2)]-\Expect_{\hat{F}}[\min_{x'}(x'-\xi)^2]}{\Expect_{\hat{F}}[(\hat{x}-\xi)^2]-\Expect_{\hat{F}}[\min_{x'}(x'-\xi)^2]}\\
%         &=\coeffDet(x(\cdot)),
%     \end{aligned}
% \end{equation*}
since $\Expect_{\hat{F}}[\min_{x'}(x'-\xi)^2]=0$ and $\hat{x}:=\argmin_{\x}\Expect_{\hat{F}}[(x-\xi)^2]=\hat{\xi}$. Hence, the coefficient of prescriptiveness has a similar interpretation as $\coeffDet$. Namely,  $\coeffPresc_{F}$ is upper bounded by $1$, and as it gets closer to $1$, it indicates how successful the data-driven policy has been in closing the gap between the SAA solution that makes no use of covariate information and a hypothetical policy that would have access to full information about $\Vxi$. 

One can also find traces in the literature of attempts to measure $\coeffDet(x(\cdot))$ out-of-sample. Namely, \cite{10.2307/40056860} studies whether excess stock return predictors can outperform historical averages in terms of out-of-sample explanatory power of such predictors. This measure can be captured using \[R_F^2(f(\cdot),\hat{F}) := 1- \frac{\Expect_{F}[(f(
\Vzeta)-\xi)^2]}{\Expect_{F}[(\hat{\xi}-\xi)^2]}=\coeffPresc_F(f(\cdot))\,,\]
which naturally leads to the question of whether $R^2(f(\cdot))$ is a good approximation for $R_F^2(f(\cdot),\hat{F})$ in a data-driven environment (potentially susceptible to distribution shifts). If not, then one must turn to employing more robust estimation methods.

\section{Robust Data-driven Prescriptiveness Optimization}\label{sec:Model}

In order to tackle the robustification and optimization of $\coeffPresc$, we consider a more general version of this measure, which relaxes the assumption that the benchmark is the solution to \eqref{eq:SP_data} and widens the scope of our analysis. To this end, we define the prescriptiveness competitive ratio (PCR) of a policy $\x(\cdot)$ with respect to a reference policy $\bar{x}$ as:
% \begin{equation}
%  \resizebox{1\hsize}{!}{$
%     \ratioPresc_F(\x(\cdot),\bar{\x}):= \left\{\begin{array}{ll}1- \frac{\Expect_F[h(\x(\Vzeta),\V{\xi})]-\Expect_F[\min_{\x'\in\X}h(\x',\V{\xi})]}{\Expect_F[h(\bar{\x},\V{\xi})]-\Expect_F[\min_{\x'\in\X}h(\x',\V{\xi})]} & \mbox{if $\Expect_F[h(\bar{\x},\V{\xi})]-\Expect_F[\min_{\x'\in\X}h(\x',\V{\xi})]>0$}\\
%     1 & \hspace{-2.2cm}\mbox{if $\Expect_F[h(\bar{\x},\V{\xi})]=\Expect_F[\min_{\x'\in\X}h(\x',\V{\xi})]=\Expect_F[h(\x(\Vzeta),\V{\xi})]$}\\
%     -\infty & \hspace{-2.2cm}\mbox{otherwise}\end{array}\right.\hspace{-0.25cm}.
% $}
% \end{equation}
\begin{equation}
\begin{aligned}
    &\ratioPresc_F(\x(\cdot),\bar{\x}):=
    &\left\{\hspace{-1ex}\begin{array}{l}1-\frac{\Expect_F[h(\x(\Vzeta),\V{\xi})]-\Expect_F[\min_{\x'\in\X}h(\x',\V{\xi})]}{\Expect_F[h(\bar{\x},\V{\xi})]-\Expect_F[\min_{\x'\in\X}h(\x',\V{\xi})]}\vspace{1ex}\\
    \qquad\;\;\mbox{if}\begin{array}{c}
    \Expect_F[h(\bar{\x},\V{\xi})]
    -\Expect_F[\min_{\x'\in\X}h(\x',\V{\xi})]>0
    \end{array}\vspace{2ex}\\
    1\quad\;\;\;\; \mbox{if }  \Expect_F[h(\bar{\x},\V{\xi})]=\Expect_F[\min_{\x'\in\X}h(\x',\V{\xi})]
    \vspace{2ex}\\
    \qquad\;\;\mbox{and } \Expect_F[h(\x(\Vzeta),\V{\xi})]=\Expect_F[\min_{\x'\in\X}h(\x',\V{\xi})]\vspace{2ex}\\
    -\infty\;\;\;\mbox{otherwise}\end{array}\right.\hspace{-0.25cm}.
\end{aligned}
\end{equation}
Indeed, the coefficient of prescriptiveness can be considered a special case when $\bar{x}:=\xSAA$:
\[\ratioPresc_F(\x(\cdot),\xSAA)=1- \frac{\Expect_F[h(\x(\Vzeta),\V{\xi})]-\Expect_F[\min_{\x'\in\X}h(\x',\V{\xi})]}{\Expect_F[h(\xSAA,\V{\xi})]-\Expect_F[\min_{\x'\in\X}h(\x',\V{\xi})]} = \coeffPresc_F(\x(\cdot)).\]
% \begin{equation*}
% \resizebox{1\hsize}{!}{$
%     \begin{aligned}
%         \ratioPresc_F(\x(\cdot),\xSAA)&=1- \frac{\Expect_F[h(\x(\Vzeta),\V{\xi})]-\Expect_F[\min_{\x'\in\X}h(\x',\V{\xi})]}{\Expect_F[h(\xSAA,\V{\xi})]-\Expect_F[\min_{\x'\in\X}h(\x',\V{\xi})]}\\ &= \coeffPresc_F(\x(\cdot)).
%     \end{aligned}
% $}
% \end{equation*}
when $\Expect_F[h(\xSAA,\Vxi)]-\Expect_F[\min_{\x'\in\X}h(\x',\Vxi)]>0$, while the two other cases follow from the natural extension of the definition of $\coeffPresc_F(\x(\cdot))$. In contrast to $\coeffPresc_F(\x(\cdot))$ which benchmarks policy $\x(\cdot)$ only to the SAA solution, the definition of $\mathcal{V}_F$ allows to benchmark against any other static policy. This allows our model to accommodate situations where more sophisticated statistical tools might be used to obtain the reference decision {\color{black}(e.g. regularized or distributionally robust SAA approaches \citep{doi:10.1287/opre.2018.1786,mohajerin2018data,doi:10.1287/mnsc.2020.3678}, variance-based regularized solution schemes \citep{DuchiDRLossesLatentCovMix}, or data-pooled solutions schemes \citep{doi:10.1287/mnsc.2020.3933}).} \footnote{In fact, one can go a step further and define $ \ratioPresc_F(\x(\cdot),\bar{\x}(\cdot))$ were $\bar{\x}(\cdot)$ is not a static policy. For example, $\bar{\x}(\cdot)$ could be a simple rule-based policy such as the order-up-to policy in inventory control. For ease of exposition, we treat the benchmark policy $\bar{\bm x}$ as a static policy for the remainder of the paper.} 

In a finite sample regime, where $\hat{F}$ might fail to capture the true underlying distribution, or in a situation where we expect distribution shifts, one should be interested in a distributionally robust estimation of the PCR (or equivalently of the coefficient of prescriptiveness), which takes the form of:
\[\ratioPresc_\Dist(\x(\cdot),\bar{\x}):= \inf_{F\in\Dist} \ratioPresc_F(\x(\cdot),\bar{\x})=\inf_{F\in\Dist}\coeffPresc_F(\x(\cdot)) \mbox{ when $\bar{\x}:=\xSAA$}.\]
% \begin{equation*}
%     \begin{aligned}
%         \ratioPresc_\Dist(\x(\cdot),\bar{\x}):= \inf_{F\in\Dist} \ratioPresc_F(\x(\cdot),&\bar{\x})=\inf_{F\in\Dist}\coeffPresc_F(\x(\cdot))\\
%         &\mbox{ when $\bar{\x}:=\xSAA$}.
%     \end{aligned}
% \end{equation*}
PCR where $\Dist$ is a set of distribution over the joint space $(\Vzeta,\Vxi)$, and the notation $\ratioPresc_\Dist$ is overloaded to denote the distributional robust PCR measure. Furthermore, one might be interested in identifying the policy that maximizes the PCR in the form of the following distributionally robust optimization problem:
\[(DRPCR)\;\;\;\max_{\x(\cdot)\in\policySet} \ratioPresc_\Dist(\x(\cdot),\bar{\x})\]
where $\policySet\subseteq \{\x:\Re^{n_\zeta}\rightarrow\X\}$. The following lemma provides interpretable bounds for the value of $\mathcal{V_\mathcal{D}}$.

\begin{lemma}\label{thm:ratioInterval}
    If $\bar{\x}\in\policySet$, then the optimal value of DRPCR is necessarily in the interval $[0,\,1]$.
\end{lemma}

% \begin{proof}
% This follows simply from $\ratioPresc_F(\x(\cdot),
% \bar{\x})$ being bound above by 1 for all policy $x(
% \cdot)$ and all distribution $F$ due to:
% \begin{align*}  
% \ratioPresc_F(\x(\cdot),
% \bar{\x}) &= 1- \frac{\Expect_F[h(\x(\Vzeta),\V{\xi})]-\Expect_F[\min_{\x'\in\X}h(\x',\V{\xi})]}{\Expect_F[h(\bar{\x},\V{\xi})]-\Expect_F[\min_{\x'\in\X}h(\x',\V{\xi})]} \\
% &\leq 1- \frac{\Expect_F[\min_{\x'\in\X}h(\x',\Vxi)]-\Expect_F[\min_{\x'\in\X}h(\x',\V{\xi})]}{\Expect_F[h(\bar{\x},\V{\xi})]-\Expect_F[\min_{\x'\in\X}h(\x',\V{\xi})]} = 1
% \end{align*}
% when $\Expect_F[h(\xSAA,\Vxi)]-\Expect_F[\min_{\x'\in\X}h(\x',\Vxi)]>0$, and otherwise equal to 1 or $-\infty$ both bounded above by 1. Hence, 
% \[\max_{\x(\cdot)\in\policySet} \inf_{F\in\Dist}\ratioPresc_F(\x(\cdot),\bar{\x}) \leq 1.\]
% Moreover, if $\bar{\x}\in\policySet$, then we have that
% \[\max_{\x(\cdot)\in\policySet} \ratioPresc_\Dist(\x(\cdot),\bar{\x})\geq \ratioPresc_\Dist(\bar{\x},\bar{\x})=\left\{\begin{array}{cl}0&\mbox{if $\Expect_F[h(\bar{\x},\V{\xi})]-\Expect_F[\min_{\x'\in\X}h(\x',\V{\xi})]>0$}\\ 1 & \mbox{ otherwise.}\end{array}\right.\]
% \end{proof}

Lemma~\ref{thm:ratioInterval} can be interpreted as follows. First, if $\bm x(\bm \zeta)$ achieves a $\ratioPresc_\Dist(\x(\cdot),
\bar{\x})=1$ then the policy is guaranteed to exploit $\Vzeta$ just as efficiently as if it had full information about $\Vxi$ (namely achieves the fully anticipative performance). On the other end of the spectrum, $\ratioPresc_\Dist(\x(\cdot),
\bar{\x})=0$ indicates that the policy can potentially fail to exploit any of the information present in $\bm \zeta$. When  $\bar{\x}\in\policySet$, one can always prevent negative PCR by falling back to the benchmark policy $\bar{\x}$.

Next, we show that in an environment where the distribution is known,  the optimal policy obtained from CSO is an optimal solution to DRPCR. Before proceeding, we first make the following assumption.
\begin{assumption}\label{ass:fullPolicySet}
    The policy set $\policySet$ contains all possible mappings, i.e. $\policySet:=\{\x:\Re^{n_\zeta}\rightarrow \X\}$.
\end{assumption}

\begin{lemma}\label{thm:ERMreduction}
    Given that Assumption~\ref{ass:fullPolicySet} is satisfied, if the distribution set is a singleton, i.e. $\Dist=\{\bar{F}\}$, then the optimal policy obtained from the CSO problem that employs $\bar{F}$ maximizes DRPCR.
\end{lemma}

While Lemma \ref{thm:ERMreduction} implies that DRPCR reduces to CSO when the distribution is known thus making the question of PCR optimization and performance irrelevant, this is not the case anymore for larger ambiguity sets $\Dist$. 

In this section, we first present a convex reformulation of DRPCR and then provide a reformulation of the problem for the nested CVaR ambiguity set. Finally, we propose a decomposition algorithm for solving the problem based on a bisection algorithm.

\subsection{Convex formulation for DRPCR}

The following proposition provides a convex reformulation of DRPCR. 
\begin{proposition}\label{thm:epigraph}
Given that $\bar{\x}\in\policySet$, DRPCR is equivalent to
\begin{subequations}\label{eq:thmEpi:mainProb}
\begin{eqnarray}    
    \max_{\x(\cdot)\in\mathcal{H},\gamma} \;\;\;&& \gamma\\
    \subto \;\;\;&& Q(\x(\cdot),\gamma)\leq 0\\
    &&0\leq \gamma\leq 1.
\end{eqnarray}
\end{subequations}
where
\[Q(\x(\cdot), \gamma):=\sup_{F\in\mathcal{D}}\;\Expect_{{F}}\Big[h(\x(\Vzeta),\V{\xi})-\Big((1-\gamma) h(\bar{\x},\V{\xi})+\gamma \min_{\bm x'\in\X}h(\x',\V{\xi})\Big)\Big]\]
% \begin{equation*}
%     \begin{aligned}
%     Q(\x(\cdot), \gamma):=\sup_{F\in\mathcal{D}}\;\Expect_{{F}}\Big[&h(\x(\Vzeta),\V{\xi})-\Big((1-\gamma) h(\bar{\x},\V{\xi})\\
%     &+\gamma \min_{\bm x'\in\X}h(\x',\V{\xi})\Big)\Big]
%     \end{aligned}
% \end{equation*}
is a convex non-decreasing function of $\gamma$. Moreover, problem \eqref{eq:thmEpi:mainProb} is a convex optimization problem when $\policySet$ is convex.
\end{proposition}

From the reformulation \eqref{eq:thmEpi:mainProb} one can draw interesting insights regarding the connection of DRPCR and risk-averse regret minimization, see \cite{multistageRegret}. For $\gamma = 1$, the problem reduces to the ex-post risk-averse regret minimization problem. In contrast, for $\gamma = 0$, one can interpret the problem as regretting the performance of the policy compared to a policy with less information. In the notation of \cite{multistageRegret}, this will lead to a risk-averse regret problem with  $\Delta=-1$.

\subsection{The nested CVaR ambiguity set $\mathcal{D}$}

In the following, we consider a discrete empirical distribution $\bar{F}$ and restrict $\mathcal{D}$ to be a nested CVaR ambiguity set. This ambiguity set is motivated by the works on nested dynamic risk measures (see \cite{RIEDEL2004}, \cite{Detlefsen2005} and \cite{shapiroNestedCVaR}) as will be explained shortly. We formalize our approach through the following assumption.

\begin{assumption}\label{ass:CVaR}
There is a discrete distribution $\bar{F}$, with $\{\Vzeta_\omega\}_{\omega\in\Omega_\zeta}$ and $\{\bm\xi_\omega\}_{\omega\in\Omega_\xi}$ as the set of distinct scenarios for $\Vzeta$ and $\Vxi$ respectively, such that the distribution set $\Dist$ takes the form of the \quoteIt{nested CVaR ambiguity set} with respect to $\Prob_{\bar{F}}$ and defined as
% \begin{equation}\label{CVaR}
%     \DistCVaR(\bar{F},\alpha):=\left\{F\in\mathcal{M}(\Omega_\zeta\times\Omega_\xi)\,\middle|\,\begin{array}{c}\Prob_{F}(\Vzeta=\Vzeta_\omega)=\Prob_{\bar{F}}(\Vzeta=\Vzeta_\omega)\;\forall \omega\in\Omega_\zeta,\\\Prob_{F}(\V{\xi}=\Vxi_{\omega'}|\Vzeta_{\omega})\leq (1/(1-\alpha))\Prob_{\bar{F}}(\Vxi=\Vxi_{\omega'}|\Vzeta_{\omega})\;\;\forall\,\omega\in\Omega_\zeta,\omega'\in\Omega_\xi\end{array}\right\}.
% \end{equation}
\begin{equation}\label{CVaR}
    \DistCVaR(\bar{F},\alpha):=\left\{F\in\mathcal{M}(\Omega_\zeta\times\Omega_\xi)\,\middle|\,\begin{array}{l}\Prob_{F}(\Vzeta=\Vzeta_\omega)=\Prob_{\bar{F}}(\Vzeta=\Vzeta_\omega)\;\forall \omega\in\Omega_\zeta,\\\Prob_{F}(\V{\xi}=\V{\xi}_{\omega'}|\Vzeta_{\omega})\leq (1/(1-\alpha))\Prob_{\bar{F}}(\V{\xi}=\V{\xi}_{\omega'}|\Vzeta_{\omega}
    )\;\;\\ \hspace{4.6cm}\forall\,\omega\in\Omega_\zeta,\omega'\in\Omega_\xi\end{array}\right\}.
\end{equation}
% \begin{equation}\label{CVaR}
% \resizebox{1\hsize}{!}{$
%     \begin{aligned}
%     &\DistCVaR(\bar{F},\alpha):=\\
%     &\left\{F\in\mathcal{M}(\Omega_\zeta\times\Omega_\xi)\,\middle|\begin{array}{l}\Prob_{F}(\Vzeta=\Vzeta_\omega)=\\
%     \qquad\Prob_{\bar{F}}(\Vzeta=\Vzeta_\omega)\;\forall \omega\in\Omega_\zeta,\vspace{2ex}\\
%     \Prob_{F}(\V{\xi}=\V{\xi}_{\omega'}|\Vzeta_{\omega})\leq\\
%     \qquad(1/(1-\alpha))\Prob_{\bar{F}}(\V{\xi}=\V{\xi}_{\omega'}|\Vzeta_{\omega}
%     )\;\;\\ \qquad\forall\,\omega\in\Omega_\zeta,\omega'\in\Omega_\xi\end{array}\hspace{-2ex}\right\}.
%     \end{aligned}
% $}
% \end{equation}
where $\mathcal{M}(\Omega_\zeta\times\Omega_\xi)$ is the set of all distributions supported on over the joint space $\{\Vzeta_\omega\}_{\omega\in\Omega_\zeta}\times\{\bm\xi_\omega\}_{\omega\in\Omega_\xi}$.
\end{assumption}

The structure of $\DistCVaR(\bar{F},\alpha)$ implies that  there is no ambiguity in the marginal distribution of the  observed random variable $\bm \zeta$. Rather, the ambiguity is solely on the unobserved random variable $\bm\xi$ and is sized using the parameter $\alpha$. The nested CVaR ambiguity set owes its name from \cite{shapiroNestedCVaR} and the fact that for any function $g(\x,\Vxi)$:
\begin{align*}\label{eq:D1}
\sup_{F\in\DistCVaR(\bar{F},\alpha)}\;&\Expect_{{F}}\left[g(\x(\Vzeta),\V{\xi})\right]\\
&=\sup_{F\in\DistCVaR(\bar{F},\alpha)}\;\sum_{\omega\in\Omega_\zeta} \sum_{\omega'\in\Omega_\xi}\Prob_F(\Vzeta=\Vzeta_\omega)\cdot \Prob_F(\Vxi=\Vxi_{\omega'}|\Vzeta=\Vzeta_\omega) g(\x(\Vzeta_\omega),\Vxi_{\omega'})\\
%&\;\;=\sum_{\omega\in\Omega_\zeta} \Prob_{\bar{F}}(\Vzeta=\Vzeta_\omega) \sup_{F_{\xi|\zeta}\in \DistCVaR(\bar{F}_{\xi|\zeta_\omega},\alpha) } \sum_{\omega'\in\Omega_\xi}  \Prob_{F_{\xi|\zeta}}(\Vxi=\Vxi_{\omega'}) g(\x(\Vzeta_\omega),\Vxi_{\omega'})\\
&=\Expect_{\bar{F}}\left[\CVaR^{\alpha}_{\bar{F}}\left(g(\x(\Vzeta),\Vxi)|\Vzeta\right)\right].
\end{align*}
%where $\bar{F}_{\xi|\zeta_\omega}$ is the conditional distribution of $\bar{F}$ given $\Vzeta_\omega$ and where we overload the notation of $\DistCVaR$ letting:
%\[\DistCVaR(\bar{F}_{\xi|\zeta},\alpha):=\{F_{\xi|\zeta}\in\mathcal{M}(\Omega_\xi):\Prob_{F_{\xi|\zeta}}(\Vxi=\Vxi_{\omega'})\leq(1/(1-\alpha))\Prob_{\bar{F}_{\xi|\zeta}}(\Vxi=\Vxi_{\omega'})\forall\,\omega'\in\Omega_\xi\}.\]
For $\alpha = 0$, the problem reduces to $\min_{\bm x(\cdot)\in\mathcal{H}}\;\Expect_{\bar{F}}[h(\x(\Vzeta),\Vxi)]$, effectively recovering the CSO policy. On the other spectrum, for $\alpha = 1$ the problem reduces to $\min_{\bm x(\cdot)\in\mathcal{H}}\;\Expect_{\bar F}[\max_{\omega:\bar{P}(\Vxi=\Vxi_\omega|\Vzeta)>0} h(\x(\Vzeta),\Vxi_\omega)]$, which implies that  for each realization of $\bm \zeta_\omega$ the decision $\bm x(\bm \zeta_\omega)$ is robust against all admissible realizations of $\bm \xi$ given $\bm \zeta_\omega$. %\EDmodified{We finally note that our result do extend to any ambiguity sets that impose that the distribution of $\zeta$ and conditional distributions of $\xi$ given $\zeta$ to lie in their respective convex uncertainty set (see Appendix \ref{sec:appgenCVaR} for more details).}

The nested CVaR representation and full policy space Assumption~\ref{ass:fullPolicySet} can be exploited to optimize $Q(\x(\cdot),\gamma)$. 

\begin{proposition}\label{thm:epigraph_reformulation}
Under Assumption \ref{ass:CVaR}, problem~\eqref{eq:thmEpi:mainProb} can thus be reformulated as
\begin{subequations}\label{nesterCVaRproblem}
\begin{eqnarray}
    \max_\gamma&&\gamma\\
    \subto&&\displaystyle \sum_{\omega\in\Omega_\zeta} \Prob_{\bar{F}}(\Vzeta=\Vzeta_\omega) \phi_\omega(\gamma)\leq 0 \label{nesterCVaRproblemB}\\
    && 0\leq \gamma\leq 1\,,
\end{eqnarray}
\end{subequations}
where $\phi_\omega(\gamma)$ is a non-decreasing function (when $\bar{\x}\in\X$) capturing the optimal value of:
\begin{subequations}\label{eq:probScenBased}
\begin{eqnarray}
    \min_{\x\in\X,t,\V{s}\geq 0}&&t+\frac{1}{1-\alpha}\sum_{\omega'\in\Omega_\xi} \Prob_{\bar{F}}(\V{\xi}=\V{\xi}_{\omega'}|\Vzeta=\Vzeta_\omega)s_{\omega'}\\
    \subto && s_{\omega'} \geq h(\x,\V{\xi}_{\omega'})-\left((1-\gamma) h(\bar{\x},\V{\xi}_{\omega'})+\gamma \min_{x'\in\X}h(\x',\V{\xi}_{\omega'})\right)- t\,,\;\forall \omega'\in\Omega_\xi\quad\quad\;\label{eq:probScenBased:C2}.
\end{eqnarray}
\end{subequations}
% \begin{subequations}\label{eq:probScenBased}
% \begin{align}
%     \min_{\x\in\X,t,\V{s}\geq 0}&\quad t+\frac{1}{1-\alpha}\sum_{\omega'\in\Omega_\xi} \Prob_{\bar{F}}(\V{\xi}=\V{\xi}_{\omega'}|\Vzeta=\Vzeta_\omega)s_{\omega'}\\
%     \subto &\quad s_{\omega'} \geq h(\x,\V{\xi}_{\omega'})-\Big((1-\gamma) h(\bar{\x},\V{\xi}_{\omega'})\notag\\
%     &\;\quad+\gamma \min_{x'\in\X}h(\x',\V{\xi}_{\omega'})\Big)- t,\,\;\forall \omega'\in\Omega_\xi \label{eq:probScenBased:C2}
% \end{align}
% \end{subequations}
and can be reduced to a linear program when $\X$ is polyhedral and $h(\x,\Vxi_{\omega'})$ is linear programming representable. 
\end{proposition}

\EDmodified{In practice, $\bar{F}$ is often composed of an empirical distribution $\hat{F}_\zeta$ and a trained conditional distribution $\hat{F}_{\xi|\zeta}$. Given an optimal solution $\gamma^*$ to problem \eqref{nesterCVaRproblem}, one should then define the extended optimal policy $\x(\Vzeta)$ beyond $\{\Vzeta_\omega\}_{\omega\in\Omega_\zeta}$ using the optimal solution of problem \eqref{eq:probScenBased} with $\hat{F}_{\xi|\zeta}$.}

\pushToAppendix{
\begin{proof}
Letting $g(\x,\Vxi,\gamma):=h(\x,\V{\xi})-\left((1-\gamma) h(\bar{\x},\V{\xi})+\gamma \min_{x'\in\X}h(\x',\V{\xi})\right)$, we have that
\begin{align*}  
\psi(\gamma):=&\min_{\x(\cdot)\in\mathcal{H}}Q(\x(\cdot),\gamma)\\
=&\min_{\x(\cdot)\in\mathcal{H}}\sup_{F\in\DistCVaR(\bar{F},\alpha)}\;\Expect_{{F}}\Big[g(\x(\Vzeta),\V{\xi},\gamma)\Big]\\
=&\min_{\x(\cdot)\in\mathcal{H}}\Expect_{\bar{F}}\Bigg[\CVaR^{\alpha}_{\bar{F}}\Big(g(\x(\Vzeta),\Vxi,\gamma)|\Vzeta\Big)\Bigg]\\
=&\min_{\x(\cdot)\in\mathcal{H}}\Expect_{\bar{F}}\Bigg[\inf_t\;t+\frac{1}{1-\alpha}\Expect_{\bar{F}}\bigg[\max\Big(0,g(\x(\Vzeta),\Vxi,\gamma)-t\Big)|\Vzeta\bigg]\Bigg]\\
=&\Expect_{\bar{F}}\Bigg[\inf_{\x\in\X,t}\;t+\frac{1}{1-\alpha}\Expect_{\bar{F}}\bigg[\max\Big(0,g(\x,\Vxi,\gamma)-t\Big)|\Vzeta\bigg]\Bigg],
\end{align*}
where we exploit the infimum representation of CVaR and the interchangeability property of expected value operators (see \cite{SHAPIRO2017377} and reference therein).
Given that $\bar{F}$ is a discrete distribution as described in Assumption \ref{ass:CVaR}, one can compute $\psi(\gamma)$ by solving for each scenario $\Vzeta_\omega$ with $\omega\in\Omega_\zeta$ the problem \eqref{eq:probScenBased:C2}.
% \begin{subequations}\label{eq:probScenBased}
% \begin{eqnarray}
%     \phi_\omega(\gamma):=\quad\quad\min_{\x\in\X,t,\V{s}}&&t+\frac{1}{1-\alpha}\sum_{\omega'\in\Omega_\xi} \Prob_{\bar{F}}(\Vxi=\Vxi_{\omega'}|\Vzeta=\Vzeta_\omega)s_{\omega'}\\
%     \subto && s_{\omega'} \geq h(\x,\V{\xi}_{\omega'})-\left((1-\gamma) h(\bar{\x},\V{\xi}_{\omega'})+\gamma \min_{x'\in\X}h(\x',\V{\xi}_{\omega'})\right) - t\,,\;\forall \omega'\in\Omega_\xi\label{eq:probScenBased:C2}\\
%     &&s_{\omega'} \geq 0\,,\;\forall \omega'\in\Omega_\xi.
% \end{eqnarray}
% \end{subequations}
Based on the solution of problem \eqref{eq:probScenBased} for each $\omega\in\Omega_\zeta$, one can obtain $\psi(\gamma):=\sum_{\omega\in\Omega_\zeta}\Prob_{\bar{F}}(\Vzeta=\Vzeta_\omega) \phi_\omega(\gamma)$ together with a potentially feasible policy $\x(\Vzeta):=\x_{\omega(\Vzeta)}^*$, 
where $\omega(\Vzeta)=\argmin_{\omega\in\Omega_\zeta}\|\Vzeta-\Vzeta_\omega\|$ and $\x_\omega$ refers to the minimizer of problem \eqref{eq:probScenBased}.
We further note that problem \eqref{eq:probScenBased} can be reduced to a linear program when $\X$ is polyhedral and $h(\x,\Vxi_{\omega'})$ is linear programming representable for all $ \omega'\in\Omega_\xi$.    
\end{proof}
}

This being said, whether problem \eqref{eq:probScenBased} is reduceable to a linear program or, more generally, a convex optimization model, its size scales with $|\Omega_\zeta|\cdot|\Omega_\xi|$, which can be computationally challenging. 
\AGmodified{
We therefore propose a decomposition algorithm to efficiently solve the problem.   Let $\psi(\gamma) := \sum_{\omega\in\Omega_\zeta} \Prob_{\bar{F}}(\Vzeta=\Vzeta_\omega) \phi_\omega(\gamma)$. Using the definition of $\phi_\omega(\gamma)$, we observe that for fixed $\gamma$ one can evaluate $\psi(\gamma)$ by solving $|\Omega_\zeta|$ distinct problem \eqref{eq:probScenBased} for each $\omega\in\Omega_\zeta$. Moreover, given that each $\phi_\omega(\gamma)$ is non-decreasing (see Proposition \ref{thm:epigraph_reformulation}), one concludes that  %Proposition~\ref{thm:epigraph} 
$\psi(\gamma)$ is non-decreasing. Hence, one can design a bisection algorithm on $\gamma$ to solve the DRPCR problem \eqref{eq:thmEpi:mainProb}.  Namely, each step consists in identifying the mid-point $\tilde{\gamma}$ of an interval known to contain the optimal value of $\gamma$, and verifying whether $\tilde{\gamma}$ is feasible by evaluating $\psi(\tilde \gamma)$ to decide which of the two sub-interval below or above $\tilde{\gamma}$ contains $\gamma^*$, see Figure~\ref{fig:bisection} (left) in Appendix~\ref{Appendix:Bisection}. The details of this algorithm are presented in Algorithm~\ref{alg:bisection}. It's efficiency relies on the difficulty of executing step \ref{step:bisAlg:solving}, i.e. evaluation $\phi_\omega(\gamma)$ for each $\omega$. 
{\color{black}
The following lemma provides formal guaranties regarding the convergence rate of Algorithm~\ref{alg:bisection}.
\begin{lemma}\label{convergence_algorithm1}
    Algorithm~\ref{alg:bisection} terminates in $\lceil\log_2(1/\epsilon)\rceil$ iterations. Moreover, if $\mathcal{X}$ is polyhedral and $h(\bm x,\bm\xi)$ linear programming representable, the algorithm terminates in polynomial time with respect to $\log(1/\epsilon)$, $|\Omega_\zeta|$, $|\Omega_\xi|$, $n_x$, $n_\xi$, the size of the LP representation of $\mathcal{X}$ and of $h(\boldsymbol x,\boldsymbol\xi).$
\end{lemma}
}

Appendix~\ref{Appendix:Bisection} further proposes an accelerated bisection algorithm for  the case when $\mathcal{X}$ is convex. Namely, it derives the sub-gradient of $\psi(\gamma)$ and exploits its convexity to tighten the interval for $\gamma^*$ at each iteration.  %{\color{red}Note that in the case where $\mathcal{X}$  is non-polyhedral but convex, typically near optimal solutions are achieved for problem \eqref{eq:probScenBased}. In such case, modifications to the bisection algorithm can be applied, see for example \cite{CohenLMPS16}.}
%\EDmodified{Note that in the case where problem \eqref{eq:probScenBased} can only be approximately solved, e.g. when $\mathcal{X}$ is convex but non-polyhedral, a slightly modified version of our our bisection algorithm might be needed, see for example \cite{CohenLMPS16}.}
}

\begin{algorithm}[H]
\caption{Bisection algorithm for DRPCR}\label{alg:bisection}
\begin{algorithmic}[1]
\State Input: Tolerance $\epsilon>0$
\State Set $\gamma^-:=0$, $\gamma^+:=1$
\While{$\gamma^+-\gamma^->\epsilon$}
\State Set $\tilde{\gamma}:=(\gamma^+ +\gamma^-)/2$
\State //Solve $\min_{\x(\cdot)\in\policySet} Q(\x(\cdot),\tilde\gamma)$ to get %optimal policy $\tilde{\x}(\cdot)$ and 
optimal value $\tilde{\psi}(\tilde{\gamma})$ 
\For{$\omega\in\Omega_\zeta$}
\State Solve problem \eqref{eq:probScenBased} with $\omega$ and $\tilde{\gamma}$ to get optimal value $\phi_{\omega}(\tilde{\gamma})$\label{step:bisAlg:solving}\label{step:solve_subproblem}
\EndFor
%\State Let $\tilde{\psi}:=\sum_{\omega\in\Omega_\zeta} \Prob_{\bar{F}}(\Vzeta=\Vzeta_\omega) \phi_\omega(\gamma)$
%\State {\color{blue}Evaluate $\phi_\omega(\tilde\gamma)$ and get the optimal $\tilde{\x}_\omega$ for all $\omega\in\Omega_\zeta$. Set $\tilde\psi= \sum_{\omega\in\Omega_\zeta} \Prob_{\bar{F}}(\Vzeta=\Vzeta_\omega) \phi_\omega(\tilde\gamma)$}\label{step:bisAlg:solving}
\If{$\tilde{\psi}:=\sum_{\omega\in\Omega_\zeta} \Prob_{\bar{F}}(\Vzeta=\Vzeta_\omega) \phi_\omega(\tilde{\gamma})\leq 0$}\label{step:prob11}
\State Set $\gamma^-:=\tilde{\gamma}$% and %{\color{blue}$\x^*(\bm \zeta_\omega):=\tilde{\x}_\omega$ for all $\omega\in\Omega_\zeta$}  
%$\x^*(\cdot)=\tilde{\x}(\cdot)$
\Else
\State Set $\gamma^+:=\tilde{\gamma}$
\EndIf
\EndWhile
\State Return $\gamma^*:=\gamma^-$ %and $\bm x^*(\cdot)$ 
\end{algorithmic}
\end{algorithm}
{\color{black}

\subsection{Generalized nested ambiguity set $\mathcal{D}$}
One can generalize the results of the previous section by considering a generalized version of the ambiguity set formalized in the following assumption. 
\begin{assumption}\label{{ass:nested_ambiguity}}
    For a discrete distribution $\bar{F}$, the  distribution set $\Dist$ takes the form of the \quoteIt{generalized nested ambiguity set} with respect to $\Prob_{\bar{F}}$ and defined as
    \begin{equation}\label{general_nested}
 \begin{array}{ll}
\DistCVaR(\bar{F},r_\zeta,r_\xi):=
\left\{F\in\mathcal{M}(\Omega_\zeta\times\Omega_\xi)\,\middle|\begin{array}{l}d_\zeta(F_\zeta,\bar{F}_\zeta)\leq r_\zeta\\
    d_\xi\big(F_{\xi|\zeta_\omega}, \bar{F}_{\xi|\zeta_\omega}\big)\leq r_\xi  \;\;\forall\,\omega\in\Omega_\zeta\end{array}\right\}.
\end{array}
\end{equation}
where $d_\zeta(F_\zeta,\bar{F}_\zeta)$ and $d_\xi(F_{\xi|\zeta_\omega},\bar{F}_{\xi|\zeta_\omega})$ are two convex divergence measures, i.e. non-negative, convex in their first argument and minimized when the two probability measures are equal, applied on marginal distribution of $F$ and the conditional distribution of $\xi$ given $\zeta$, respectively. 
\end{assumption}
The structure of \eqref{general_nested} allows one to control the ambiguity about both the marginal distribution of $\bm\zeta$ and the conditional distributions of $\bm\xi$ using the parameters $r_\zeta$ and $r_\xi$ to bound the maximum divergence respectively. In particular, it reduces to the nested CVaR ambiguity set when using 
\begin{align*}
d_\zeta(F_\zeta,\bar{F}_\zeta):=
\inf\left\{s|\Prob_{F_\zeta}(\V{\zeta}=\V{\zeta}_\omega)\le\frac{1}{1-s}\Prob_{\bar{F}_\zeta}(\V{\zeta}=\V{\zeta}_\omega),\forall\omega\in\Omega\right\},\\
d_\xi(F_{\xi|\zeta},\bar{F}_{\xi|\zeta}):=
\inf\left\{s|\Prob_{F_{\xi|\zeta}}(\V{\xi}=\V{\xi}_{\omega'})\le\frac{1}{1-s}\Prob_{\bar{F}_{\xi|\zeta}}(\V{\xi}=\V{\xi}_{\omega'}),\forall\omega'\in\Omega'\right\},
\end{align*}
$r_\zeta:=0$, and $r_\xi:=\alpha$.

In the following, to simplify presentation, given that $\Omega_\zeta$ and $\Omega_\xi$ are finite, we let $\V{p}\in\Re^{|\Omega_\zeta|}$ denote the vector of probabilities $p_\omega:=\Prob_F(\zeta=\zeta_\omega)$ and  $\V{q}^\omega\in \Re^{|\Omega_\xi|}$ denote the probabilities $q_{\omega'}^\omega:=\Prob_F(\xi=\xi_{\omega'}|\zeta=\zeta_\omega)$, and similarly for $\bar{p}$ and $\bar{q}^\omega$ to captures the same probabilities under $\bar{F}$. We will further abuse notation and denote $d_\zeta(\V{p},\bar{\V{p}}):=d_\zeta(F_\zeta,\bar{F}_\zeta)$ and $d_\xi(\V{q}^\omega,\bar{\V{q}}^\omega):=d_\xi(F_{\xi|\zeta_\omega},\bar{F}_{\xi|\zeta_\omega})$. The following proposition generalizes Proposition~\ref{thm:epigraph_reformulation}.

\begin{proposition}\label{generalized_proposition2}
Under Assumption~\ref{{ass:nested_ambiguity}}, problem~\eqref{eq:thmEpi:mainProb} can thus be reformulated as
\begin{subequations}\label{nesterCVaRproblem}
\begin{eqnarray}
    \max_\gamma&&\gamma\\
    \subto&& \sup_{\V{p}\in\mathcal{Z}} \sum_{\omega\in\Omega_\zeta} p_\omega \bar{\phi}_\omega(\gamma) \leq 0 \label{nesterCVaRproblemB}\\
    && 0\leq \gamma\leq 1\,,
\end{eqnarray}
\end{subequations}
% \begin{equation}
% \begin{alighed}{cl}
%     \displaystyle\max_\gamma&\gamma\\
%     \subto&\displaystyle \sup_{\V{p}\in\mathcal{Z}} \sum_{\omega\in\Omega_\zeta} p_\omega \bar{\phi}_\omega(\gamma) \leq 0 \label{nesterCVaRproblemB}\\
%     & 0\leq \gamma\leq 1\,,
% \end{array}
% \end{equation}
where $\mathcal{Z}:=\{\V{p}:\V{p}\geq 0,\bm e^\top \bm p =1, d_\zeta(\V{p},\bar{\V{p}})\leq r_\zeta\}$ and $\bar{\phi}_\omega(\gamma)$ is a non-decreasing function (when $\bar{\x}\in\X$) capturing the optimal value of:
\begin{subequations}\label{eq:prop2_subproblem}
    \begin{eqnarray}
    \min_{\x,t,\alpha,\V{s}}&& t+r_\xi\alpha + d_*(\V{s},\alpha,\bar{\V{q}}^\omega)\\
    \subto && s_{\omega'} \geq h(\x,\V{\xi}_{\omega'})-\Big((1-\gamma) h(\bar{\x},\V{\xi}_{\omega'})
    +\gamma \min_{x'\in\X}h(\x',\V{\xi}_{\omega'})\Big)- t,\forall \omega'\in\Omega_\xi\quad\quad\\
    &&\x\in\X,\,\alpha\geq 0,
\end{eqnarray}
\end{subequations}
where $d_*(\bm s,\alpha,\bar{\V{q}}^\omega):=\sup_{\V{q}} \bm s^T\V{q}-\alpha d_\xi(\V{q},\bar{\V{q}}^\omega)$ is the perspective of the convex conjugate  of $d_\xi(\V{q},\bar{\V{q}}^\omega)$.
\end{proposition}
Algorithm~\ref{alg:bisection} can  be applied in the generalized setting with the simple modification that problem~\eqref{eq:probScenBased} in step~\ref{step:solve_subproblem} is replaced with the convex problem~\eqref{eq:prop2_subproblem}, and step \ref{step:prob11} must compute $\sup_{\V{p}\in\mathcal{Z}} \sum_{\omega\in\Omega_\zeta} p_\omega \bar{\phi}_\omega(\gamma)$, which now requires solving a convex optimization problem.
}

\pushToAppendix{
\subsection{A bisection algorithm for DRPCR}

From the definition of $\psi(\gamma)$ and $\phi_{\omega}(\gamma)$, we observe that for fixed $\gamma$ one can evaluate $\psi(\gamma)$ by solving $|\Omega_\zeta|$ distinct problem \eqref{eq:probScenBased} for each $\omega\in\Omega_\zeta$. Moreover, Proposition~\ref{thm:epigraph} states that $\psi(\gamma)$ is an non-decreasing convex function of $\gamma$. Hence, one can design a bisection algorithm on $\gamma$ to solve the DRPCR problem \eqref{eq:thmEpi:mainProb}.  Namely, each step consists in identifying the mid-point $\tilde{\gamma}$ of an interval known to contain the optimal value of $\gamma$, and verifying whether $\tilde{\gamma}$ is feasible by solving $\min_{\x(\cdot)\in\policySet} Q(\x(\cdot),\tilde\gamma)$ to decide which of the two sub-interval below or above $\tilde{\gamma}$ contains $\gamma^*$, see Figure~\ref{fig:bisection} (left). The details of this algorithm are presented in Algorithm \ref{alg:bisection}. It's efficiency relies on the difficulty of executing step \ref{step:bisAlg:solving}, i.e. solving $\min_{\x(\cdot)\in\policySet} Q(\x(\cdot),\gamma)$.

\begin{algorithm}[H]
\caption{Bisection algorithm for DRPCR}\label{alg:bisection}
\begin{algorithmic}[1]
\State Input: Tolerance $\epsilon>0$
\State Set $\gamma^-:=0$, $\gamma^+:=1$, $\x^*(\Vzeta):=\bar{\x}$ for all $\Vzeta$
\While{$\gamma^+-\gamma^->\epsilon$}
\State Set $\tilde{\gamma}:=(\gamma^+ +\gamma^-)/2$
\State Solve $\min_{\x(\cdot)\in\policySet} Q(\x(\cdot),\gamma)$ to get optimal policy $\tilde{\x}(\cdot)$ and optimal value $\tilde{\psi}$ \label{step:bisAlg:solving}
\If{$\tilde{\psi}\leq 0$}
\State Set $\gamma^-:=\tilde{\gamma}$ and $\x^*(\cdot):=\tilde{\x}(\cdot)$
\Else
\State Set $\gamma^+:=\tilde{\gamma}$
\EndIf
\EndWhile
\State Return $\gamma^*:=\gamma^-$ and $x^*(\cdot)$
\end{algorithmic}
\end{algorithm}

One can possibly accelerate the convergence rate on the bisection algorithm by exploiting the fact that $\psi(\cdot)$ is a non-decreasing convex function when $\X$ is convex. Indeed, for the current interval $[\gamma^-,\gamma^+]$, $\psi(\gamma)$ can be under- and over-estimated, see Figure~\ref{fig:bisection} (right). The procedure can be described as follows. First, we construct a line that will underestimate $\psi$   by identifying a subgradient of the function at $\tilde\gamma$. This can be computed analytically since 
\begin{equation*}
\begin{array}{rl}
\psi(\gamma)&:=\displaystyle\Expect_{\bar F}\left[\min_{x\in\mathcal{X}} \CVaR_\alpha\left(h(\x,\V{\xi})-\Big((1-\gamma) h(\bar{\x},\V{\xi})+\gamma \min_{x'\in\X}h(\x',\V{\xi})\Big)\middle|\zeta\right)\right]\\[2ex]

&=\displaystyle\Expect_{\bar F}\left[\min_{x\in\mathcal{X}} \sup_{F\in\DistCVaR(\bar{F},\alpha)}\Expect_F\left[h(\x,\V{\xi})-\Big((1-\gamma) h(\bar{\x},\V{\xi})+\gamma \min_{x'\in\X}h(\x',\V{\xi})\Big)\middle|\zeta\right]\right]\\[2ex]

&\geq\displaystyle\Expect_{\bar F}\left[ \sup_{F\in\DistCVaR(\bar{F},\alpha)}\min_{x\in\mathcal{X}}\Expect_F\left[h(\x,\V{\xi})-\Big((1-\gamma) h(\bar{\x},\V{\xi})+\gamma \min_{x'\in\X}h(\x',\V{\xi})\Big)\middle|\zeta\right]\right]\\   
  
&\geq\displaystyle\Expect_{\bar F}\left[ \min_{x\in\mathcal{X}}\Expect_{ F^*_{\xi|\zeta}}\left[h(\x,\V{\xi})-\Big((1-\gamma) h(\bar{\x},\V{\xi})+\gamma \min_{x'\in\X}h(\x',\V{\xi})\Big)\middle|\Vzeta\right]\right]\\   

&=\displaystyle\underbrace{\Expect_{\bar F}\left[ \min_{x\in\mathcal{X}}\Expect_{ F^*_{\xi|\zeta}}\left[h(\x,\V{\xi})-  h(\bar{\x},\V{\xi}) \right]\right]}_{\mymbox{offset \quoteIt{$a$}}} + 
\gamma\underbrace{\Expect_{\bar F}\left[\Expect_{ F^*_{\xi|\zeta}}\left[ h(\bar{\x},\V{\xi}) - \min_{x'\in\X}h(\x',\V{\xi}) \right]\right]}_{\mymbox{slope \quoteIt{$b$}}},   
\end{array}
\end{equation*}
where $F_{\xi|\zeta}^*$ is the conditional probability given $\Vzeta$ of any member (hopefully a maximizer) of $\DistCVaR(\bar{F},\alpha)$. Note that the first inequality is tight  based on 
Sion's minimax theorem (see \cite{sion58:minimax}) given that $\DistCVaR(\bar{F},\alpha)$ is compact, while the second is tight as long as $F_{\xi|\zeta}^*$ achieves the supremum. Such a maximizer can be identified using:
% \[F^*_{\xi|\zeta}\in\argmax_{F_{\xi|\zeta}\in\mathcal{M}(\Omega_\zeta):\Prob_{F_{\xi|\zeta}}(\V{\xi})\leq (1-\alpha)^{-1}\Prob_{\bar{F}}(\Vxi|\Vzeta),\,\forall\Vxi}\Expect_{F_{\xi|\zeta}
% }\left[h(\x^*_\gamma(\bm\zeta),\V{\xi})-\Big((1-\gamma) h(\bar{\x},\V{\xi})+\gamma \min_{x'\in\X}h(\x',\V{\xi})\Big)\right]\]
\[F^*_{\xi|\zeta}\in\argmax_{\tiny{
\begin{array}{cc}
F_{\xi|\zeta}\in\mathcal{M}(\Omega_\zeta):\\
\Prob_{F_{\xi|\zeta}}(\V{\xi})\leq(1-\alpha)^{-1}\Prob_{\bar{F}}(\V{\xi}|\Vzeta),\,\forall\V{\xi}
\end{array}}}\Expect_{F_{\xi|\zeta}
}\left[h(\x^*_\gamma(\bm\zeta),\V{\xi})-\Big((1-\gamma) h(\bar{\x},\V{\xi})+\gamma \min_{x'\in\X}h(\x',\V{\xi})\Big)\right]\]
where $\x^*_\gamma(\bm\zeta)$ is the minimizer of \eqref{eq:probScenBased} with $\Vzeta_\omega=\Vzeta$ since $(\bm x^*_\gamma(\cdot),F^*)$, with $F^*$ as the composition of $\bar{F}$ marginalized on $\Vzeta$ and $F^*_{\Vxi|\Vzeta}$,\footnote{Namely, $P_{F^*}(\Vxi)=P_{\bar{F}}(\Vxi)$ and $P_{F^*}(\Vxi|\Vzeta)=P_{F^*_{\xi|\zeta}}(\Vxi)$ for all $\Vzeta$.} is a saddle point of:
\[g(\x(\cdot),F):=\Expect_{{F}}\Big[h(\x(\Vzeta),\V{\xi})-\Big((1-\gamma) h(\bar{\x},\V{\xi})+\gamma \min_{x'\in\X}h(\x',\V{\xi})\Big)\Big].\]
Such a $F^*_{\Vxi|\Vzeta}$ can be obtained as a side product of solving problem \eqref{eq:probScenBased} using the optimal dual variables associated with constraint \eqref{eq:probScenBased:C2}. If we denote by $\gamma_u := a/b$ then the right bound of the interval can be updated to $\gamma^{+'} := \min(\gamma^+,\gamma_u)$. 

The second step is to construct an overestimator. If $\psi(\tilde\gamma) > 0$, then we evaluate $\psi(\gamma^{-})$ and construct the line that passes through $(\gamma^{-},\psi(\gamma^{-}))$ and $(\tilde\gamma,\psi(\tilde\gamma))$. If $\psi(\tilde\gamma) < 0$ then  we evaluate $\psi(\gamma^{+})$ and construct the line that passes through $(\gamma^{+},\psi(\gamma^{+}))$ and $(\tilde\gamma,\psi(\tilde\gamma))$. We denote the point for which the line evaluates to zero as $\gamma_o$, and update the  left bound of the interval   to $\gamma^{-'} := \max(\gamma^-,\gamma_o)$. Hence, the new interval is given by $[\gamma^{-'},\gamma^{+'}]\subseteq[\gamma^{-},\gamma^{+}]$, which would potentially significantly reducing the search space.

We conclude this section by commenting that the accelerated bisection algorithm could require up two evaluations of the $\psi$ function at each iteration instead of a single one as described in the original algorithm.

\begin{figure}[h!]
\centerline{\includegraphics[scale=0.55]{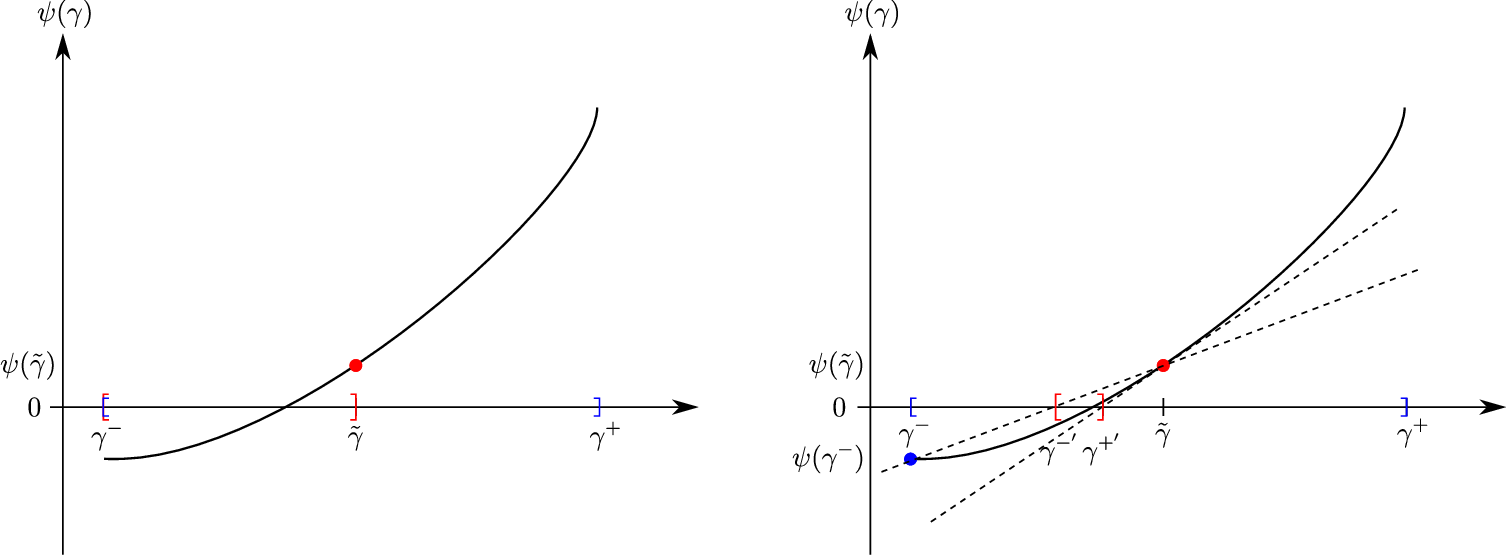}}
\caption{Visualization of the basic (left) and accelerated (right) bisection algorithm. The blue squared brackets indicate the current estimated interval containing the optimal $\gamma^*$ and the red squared brackets indicate  the interval in the next iterations. The right graph also visualizes the over and under estimators of $\psi(\gamma)$. 
}
\label{fig:bisection}
\end{figure}
}

\section{Experiments}\label{sec:Experiments}

In this section, we present a numerical study that compares the performance of DRPCR against three other data-driven benchmark methods to evaluate its robustness to perturbations of the data generating process. %More specifically, we compare the performance of the corresponding data-driven policies in terms of the coefficient of prescriptiveness over an out-of-sample dataset. 
Specifically, we will observe how these models react to the situation where one faces a distribution shift for $\V{\xi}$. In a vehicle routing problem with travel time uncertainties, this can be interpreted as a shift in the distribution of the travel times, for instance, when a special event is happening in the town. Alternatively, one can think of an inventory management problem where the manager faces a shift in the demand distribution, e.g., an unforeseen increase in demand for sanitizer during the first days of an epidemic. In general, there are numerous reasons why distribution shifts considerations might be needed depending on the context. In this regard, %this may contain many cases where one may confront so-called \quoteIt{disruptions} or \quoteIt{extreme cases} in supply chain management as an immediate result of unexpected changes in distributions of the uncertain parameters.       {\color{red} 
we refer the reader to \cite{schrouff2022maintaining} and \cite{filos2020can} for such considerations in healthcare and  autonomous driving applications.

The application that we consider for our numerical experiments is a shortest path problem described in \cite{StochasticOptimizationForests}. A directed graph is defined as $\mathcal{G}=(\mathcal{V},\,\mathcal{A})$, where $\mathcal{V}$ denotes the set of nodes and $\mathcal{A}\in \mathcal{V}\times\mathcal{V}$ is the set of arcs, i.e., ordered pairs $(i,j)$ of nodes describing the existence of a directed path from node $i$ to node $j$. The corresponding travel time of such an arc is assumed to be $\xi_{(i,j)}$. The objective of this problem is to identify the shortest path from an origin (node $o$) to a destination (node $d$). Moving away from an ideal world of known parameters gives rise to a stochastic version of this problem. In this setting, the traveling times along the arcs $\V{\xi}\in\Re^{|\mathcal{A}|}$ are uncertain; however, one might still have access to side information or observed covariates. In this case, aiming at minimizing the expected travel time leads to the following CSO problem:
\begin{equation}\label{ShortestPath}
    \x^*(\Vzeta)\in\arg\min_{\x\in\X}\Expect_{\hat{F}_{\xi|\zeta}}[\x^\top\V{\xi}],
\end{equation}
where
\[\X=\left\{\x\in\Re^{|\mathcal{A}|}\,\middle|\,
\begin{array}{l}
x_{(i,j)}\geq 0\\
\sum_{j:(o,j)\in\mathcal{A}} x_{(o,j)}-\sum_{j:(j,o)\in\mathcal{A}} x_{(j,o)}=1\\
\sum_{j:(d,j)\in\mathcal{A}} x_{(d,j)}-\sum_{j:(j,d)\in\mathcal{A}} x_{(j,d)}=-1\\
\sum_{j:(i,j)\in\mathcal{A}} x_{(i,j)}-\sum_{j:(j,i)\in\mathcal{A}} x_{(j,i)}=0
\end{array}
\begin{array}{l}
\forall (i,j)\in\mathcal{A}\\
\\%\mbox{if}\,i = o\\
\\%\mbox{if}\,i = d\\
\forall i\in\mathcal{V}\,\backslash\,\{o,\,d\}
\end{array}
\right\}\,,\]
% \begin{equation*}
% \resizebox{1\hsize}{!}{$
%     \begin{aligned}
%         \X=\left\{\x\in\Re^{|\mathcal{A}|}\,\middle|\,
%         \begin{array}{l}
%         x_{(i,j)}{\color{black}\in\{0,1\}}\remove{\geq 0}\hspace{15.5ex}\forall (i,j)\in\mathcal{A}\vspace{1ex}\\
%         \sum_{j:(o,j)\in\mathcal{A}} x_{(o,j)}-\\ \qquad\sum_{j:(j,o)\in\mathcal{A}} x_{(j,o)}=1\vspace{1ex}\\
%         \sum_{j:(d,j)\in\mathcal{A}} x_{(d,j)}-\\ \qquad\sum_{j:(j,d)\in\mathcal{A}} x_{(j,d)}=-1\vspace{1ex}\\
%         \sum_{j:(i,j)\in\mathcal{A}} x_{(i,j)}-\\ \qquad\sum_{j:(j,i)\in\mathcal{A}} x_{(j,i)}=0\,\,\,\forall i\in\mathcal{V}\,\backslash\,\{o,\,d\}
%         \end{array}
%         \hspace{-1ex}\right\}\,,
%     \end{aligned}   
% $}
% \end{equation*}
and $x_{(i,j)}=1$ if we decide to travel from node $i$ to node $j$ and $x_{(i,j)}=0$ otherwise. \remove{Similar to \cite{StochasticOptimizationForests}, we do not force integrality constraints} {\color{black}Unlike \cite{StochasticOptimizationForests}, we enforce the integrality constraints}. Furthermore, $\hat{F}_{\xi|\zeta}$ denotes the conditional distribution inferred from the training dataset.

As discussed in Section \ref{sec:intro}, DRCSO is a method proposed for robustifying the policies against distributional uncertainties in the data-driven context. Consequently, one can consider DRCSO, as an alternative to CSO, for solving this shortest-path problem. Using  the nested CVaR ambiguity set introduced in Assumption \ref{ass:CVaR} as the ambiguity set of DRCSO, one gets the model below:  
% \begin{equation}\label{eq:CVaRloss}
%     (DRCSO)\quad\quad\x^*(\Vzeta)\in\arg\min_{\x\in\X}\sup_{F_{\xi|\zeta}\in{\DistCVaR(\hat{F}_{\xi|\zeta},\alpha)}}\;\;\Expect_{F_{\xi|\zeta}}[\x^\top\V{\xi}],
% \end{equation}
\begin{equation}\label{eq:CVaRloss}
    (DRCSO)\quad\x^*(\Vzeta)\in\arg\min_{\x\in\X}\sup_{F_{\xi|\zeta}\in{\DistCVaR(\hat{F}_{\xi|\zeta},\alpha)}}\;\;\Expect_{F_{\xi|\zeta}}[\x^\top\V{\xi}],
\end{equation}
where 
% \[\DistCVaR(\hat{F}_{\xi|\zeta},\alpha):=\{F_{\xi|\zeta}\in\mathcal{M}(\Omega_\xi):\Prob_{F_{\xi|\zeta}}(\Vxi=\Vxi_{\omega'})\leq(1/(1-\alpha))\Prob_{\hat{F}_{\xi|\zeta}}(\Vxi=\Vxi_{\omega'})\forall\,\omega'\in\Omega_\xi\},\]
\begin{equation*}
    \begin{aligned}
        \DistCVaR(\hat{F}_{\xi|\zeta},\alpha):=\{F_{\xi|\zeta}\in\mathcal{M}(\Omega_\xi):\Prob_{F_{\xi|\zeta}}(\Vxi=\Vxi_{\omega'})\leq
        (1/(1-\alpha))\Prob_{\hat{F}_{\xi|\zeta}}(\Vxi=\Vxi_{\omega'})\forall\,\omega'\in\Omega_\xi\},        
    \end{aligned}
\end{equation*}
and $\alpha$ is the control parameter for the size of the ambiguity set. Staying in the DRCSO context, one can exploit a worst-case regret minimization approach instead of worst-case expected travel time. In our experiments, we look into the optimal solutions arising from an ex-post regret minimization setting, introduced as a $\Delta=1$ regret minimization model in \cite{multistageRegret}. This leads to the following distributionally robust {\color{black}contextual} regret optimization (DRCRO) problem: 
\begin{equation}\label{eq:CVaRregret}
(DRCRO)\,\,\,\x^*(\Vzeta)\in\arg\min_{\x\in\X}\sup_{F_{\xi|\zeta}\in{\DistCVaR(\hat{F}_{\xi|\zeta},\alpha)}}\;\;\Expect_{F_{\xi|\zeta}}[\x^\top\V{\xi}-\min_{x'\in\X}\x'^\top\V{\xi}].
\end{equation}
% \begin{equation}\label{eq:CVaRregret}
% \begin{aligned}
% &(DRCRO)\\
% &\x^*(\Vzeta)\in\arg\min_{\x\in\X}\sup_{F_{\xi|\zeta}\in{\DistCVaR(\hat{F}_{\xi|\zeta},\alpha)}}\;\;\Expect_{F_{\xi|\zeta}}[\x^\top\V{\xi}-\min_{x'\in\X}\x'^\top\V{\xi}].
% \end{aligned}
% \end{equation}
In this case, the decision maker compares her travel time to the one resulting from a benchmark decision that knows the future realization of $\V{\xi}$. The ultimate goal is to minimize the worst-case expectation of this gap, so-called \quoteIt{worst-case expected regret}, where the ambiguity set is nested CVaR. Finally, we solve our introduced DRPCR problem under nested CVaR ambiguity set, where the $Q(\x(\cdot),\gamma)$ function takes the form of:
\begin{equation}
    Q(\x(\cdot),\gamma):= \sup_{F\in{\DistCVaR(\hatFC,\alpha)}}\;\Expect_{{F}}\Big[\x(\Vzeta)^\top\V{\xi}-\Big((1-\gamma) \hat{\x}^\top\V{\xi}+\gamma \min_{x'\in\X}\x'^\top\V{\xi}\Big)\Big],
\end{equation}
% \begin{equation}
% \begin{aligned}
% &Q(\x(\cdot),\gamma):=\\
% &\sup_{F\in{\DistCVaR(\hatFC,\alpha)}}\;\Expect_{{F}}\Big[\x(\Vzeta)^\top\V{\xi}-\Big((1-\gamma) \hat{\x}^\top\V{\xi}+\gamma \min_{x'\in\X}\x'^\top\V{\xi}\Big)\Big],
% \end{aligned}
% \end{equation}
where $\hatFC$ denotes the distribution derived from the training dataset, composed of the empirical distribution $\hat{F}_{\zeta}$ of $\zeta$ and the inferred conditional distribution $\hat{F}_{\xi|\zeta}$,
while $\hat{\x}:= \argmin_{\x}\Expect_{\hat{F}}[h(\x,\V{\xi})]$ with $\hat{F}$ that puts equal weights on each observed data point $\{\V{\xi_i}\}_{i=1}^N$ (i.e. the SAA solution). Based on an optimal solution $\gamma^*$ for the DRPCR problem, one can retrieve an optimal policy using:
\begin{equation}
\x^*(\Vzeta)\in\arg\min_{\x\in\X}\sup_{F_{\xi|\zeta}\in{\DistCVaR(\hat{F}_{\xi|\zeta},\alpha)}}\;\;\Expect_{F_{\xi|\zeta}}\left[\x^\top\V{\xi}-\Big((1-\gamma^*) \hat{\x}^\top\V{\xi}+\gamma^* \min_{x'\in\X}\x'^\top\V{\xi}\Big)\right],\label{eq:DRPCRpolicy}
\end{equation}
% \begin{equation}
% \begin{aligned}
% \x^*(\Vzeta)\in\arg\min_{\x\in\X}\sup_{F_{\xi|\zeta}\in{\DistCVaR(\hat{F}_{\xi|\zeta},\alpha)}}\;\;\Expect_{F_{\xi|\zeta}}\bigg[\x^\top\V{\xi}-\\
% \Big((1-\gamma^*) \hat{\x}^\top\V{\xi}+\gamma^* \min_{x'\in\X}\x'^\top\V{\xi}\Big)\bigg],\label{eq:DRPCRpolicy}
% \end{aligned}
% \end{equation}
which can be obtained by solving \eqref{eq:probScenBased} with $\gamma^*$ and replacing $\Prob_{\bar{F}}(\Vxi=\Vxi_{\omega'}|\Vzeta=\Vzeta_\omega)$ with $\Prob_{\hat{F}_{\xi|\zeta}}(\Vxi_{\omega'})$. 

We adapt our numerical experiments to the graph ($\mathcal{G}$) structure employed in \cite{StochasticOptimizationForests} with the same origin ($o$) and destination ($d$); therefore, we study a graph with the size of 45 nodes ($|\mathcal{V}|=45$) and 97 arcs ($|\mathcal{A}|=97$). We assume there exist 200 covariates ($n_\Vzeta=200$) and the vector composed of travel times $\Vxi$ and covariates $\Vzeta$ follow a multivariate normal distribution. Specifically, each covariate $\zeta_i$ follows a normal distribution with a mean of zero and standard deviation of one (i.e. $\zeta_i\sim \mathcal{N}(0,\,1)$). Similarly, each travel time $\xi_{(i,j)}$ is normal with a standard deviation that matches the deviation present in \cite{StochasticOptimizationForests}'s dataset yet both the correlation and mean vector are treated differently. Starting with correlation, we introduce a new correlation structure for $(\Vzeta,\Vxi)$\footnote{This was done after observing that with \cite{StochasticOptimizationForests}'s dataset the optimal uninformed decisions produced nearly the same performance as the optimal hindsight decisions that exploited full information about realized travel times.} by instantiating a random correlation matrix (see Appendix \ref{CovMat} for details). 

Our treatment of the mean of $\Vxi$ embodies our objective to study robustness to distribution shifts. Namely, while the data generating process for the training set employs the same mean vector as in \cite{StochasticOptimizationForests}, our validation data set and out-of-sample test set will measure the performance of proposed policies on generating processes where the mean of $\Vxi$ as been perturbed, i.e. $\Expect[\xi_{(i,j)}]:=(1+\delta_{(i,j)})\mu_{(i,j)}$. Six tests were conducted for different levels of mean perturbations: no distribution shift $\delta_{(i,j)}=0$, which does not allow for any perturbation, along with tests that take into account shifts with $\delta_{(i,j)}$ generated i.i.d. according to a uniform distribution on $[0\%,\, m]$, where $m\in\mathcal{M}:=\{20\%,\,30\%,\,40\%,\,50\%,\,60\%\}$ represents the maximum possible perturbation. Furthermore, the perturbation experienced in the validation set is independent of the test set. This is to simulate situations where the level of robustness would be calibrated on a data set where a distribution shift of similar size is observed as the shift experienced out-of-sample.

\begin{figure}[htbp!]
\centerline{\includegraphics[width=\textwidth]{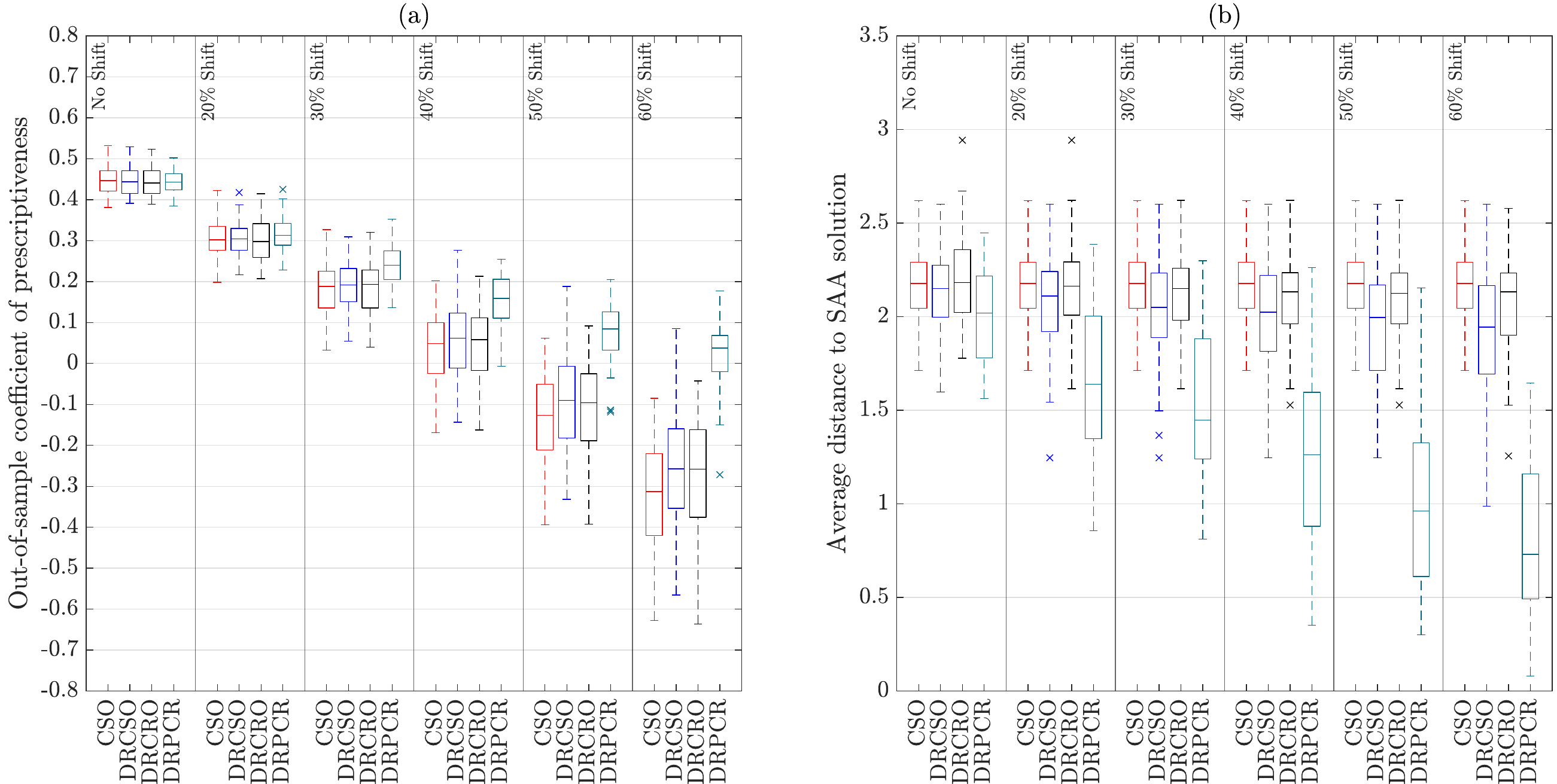}}
\caption{Shortest path problem: (a) statistics of the out-of-sample coefficient of prescriptiveness {\color{black} (lower values indicate worse performance).} (b) statistics of $\mathbb{E}_{\EDmodified{\breve{F}}}[\|\boldsymbol{x^*}(\boldsymbol \zeta) - \boldsymbol{\hat x}\|_1]$ where \EDmodified{$\breve F$ is the} out-of-sample distribution {\color{black} (lower values reflect a closer proximity to the SAA solution).}}
\label{fig:OOS}
\end{figure}

Experiments for each perturbation range contain 50 instances generated by resampling the training, validation, and test data sets. Both the training and validation datasets consist of 400 data points, while the test set contains 1000 data points and is used to measure the \quoteIt{out-of-sample} performance. The training dataset is used for learning purposes, which allows us to infer the conditional probabilities of $\hat{F}_{\xi|\zeta}$ once a new covariate vector $\Vzeta$ is observed. From a wide range of existing predictive tools for inference of $\hat{F}_{\xi|\zeta}$, \cite{BertsimasKallus} compare methods such as k-nearest-neighbors regression (\citealt{trevor2001friedman}), local linear regressions (\citealt{cleveland1988locally}), classification and regression trees (CART; \citealt{breiman1984classification}), and random forests (RF; \citealt{breiman2001random}). In their experiments, the best coefficient of prescriptiveness belongs to random forests. We exploit the code provided in \cite{StochasticOptimizationForests} to train random forests over our training datasets and then use it as the conditional distribution estimator $\hat{F}_{\xi|\zeta}$ for our validation and out-of-sample data points. The validation dataset is used to calibrate the size of the ambiguity set ($\alpha$) for the DRCSO, DRCRO, and DRPCR models. The procedure for calibrating $\alpha$ and the associated optimal $\gamma$ for the DRPCR model and to calibrate $\alpha$ for the DRCSO and DRCRO models are described in Appendix \ref{Appendix:Algo} (see algorithms \ref{alg:alpha-gamma-Cal} and \ref{alg:alpha-Cal} respectively). %Algorithm \ref{alg:alpha-gamma-Cal}. 
%A similar procedure is used to calibrate $\alpha$ for the DRCSO and DRCRO models (see Algorithm \ref{alg:alpha-Cal} \AGmodified{in Appendix~\ref{Appendix:Algo}}). 
We define the set of discretized $\alpha$ values as $\mathcal{A}:=\mathcal{A}_1\cap\mathcal{A}_2$, where $\mathcal{A}_1$ includes 20 logarithmically spaced values in $[0.01,\,0.99]$ and $\mathcal{A}_2$ includes 20 evenly spaced values in $[0,\,1)$. For CSO, Algorithm \ref{alg:alpha-Cal} can also be used with $\mathcal{A} = \{0\}.$ From a computational point of view, the training of the DRPCR algorithm took on average less than 36 minutes per instance, compared to closer to 3 minutes for DRCSO and DRCRO. The main difference comes from the extra Step \ref{step:DRPCRsolve} of Algorithm \ref{alg:alpha-gamma-Cal}, which requires solving the DRPCR problem for each candidate for $\alpha$, and took on average 50 seconds to solve, and needs to be repeated for all $\alpha\in\mathcal{A}$. {\color{black} Once the optimal $\alpha^*$ and $\gamma^*$ are determined by Algorithm \ref{alg:alpha-gamma-Cal} for a given training/validation dataset, equation \eqref{eq:DRPCRpolicy} provides the optimal policy $x^*(\cdot)$ for any covariate $\Vzeta$ received in real-time. Similarly, equations \eqref{eq:CVaRloss} and \eqref{eq:CVaRregret} can be employed to derive the optimal DRCSO and DRCRO policies associated with the real-time input of $\Vzeta$, relying on the calibrated values of $\alpha^*$ obtained from Algorithm \ref{alg:alpha-Cal}.} All optimization problems are implemented in Python and solved using Gurobi 8.1.1 on a machine featuring an Intel processor Xeon(R) CPU E5-2687W v3 @ 3.10GHz 3.10 GHz (2 processors) and 128 GB RAM. The code used for the numerical experiments is available at \href{https://github.com/erickdelage/robust_prescriptive_opt}{https://github.com/erickdelage/robust\_prescriptive\_opt}.%\href{https://anonymous.4open.science/r/Robust-Data-driven-Prescriptiveness-Optimization-9F6F}{https://anonymous.4open.science/r/Robust-Data-driven-Prescriptiveness-Optimization-9F6F}.

%\MPcomments{A main step in calibrating the size of the ambiguity set $(\alpha)$ concerning DRPCR, DRCSO, and DRCRO methods consists in the validation phase, described as Step 7 of Algorithm \ref{alg:alpha-gamma-Cal} and Step 4 of Algorithm \ref{alg:alpha-Cal}. The average runtime of this step across all instances and perturbation levels remains less than 3 minutes, for any of these three methods regardless of a binary or relaxed setting; however, one should note that the additional price of obtaining DRPCR policies, compared to the alternative ones, mainly consists in the computations embedded in Step 6 of Algorithm \ref{alg:alpha-gamma-Cal}. The average runtime of this step across all instances and perturbation levels is 27 minutes for the relaxed version of the experiments and 33 minutes for the non-relaxed problem. One should note that these cumulative runtimes reflect the total time taken for each step to be repeated across all exploited $\alpha$ levels.} 

\pushToAppendix{
\begin{algorithm}[htbp]
\caption{Algorithm of calibrating the size of the ambiguity set ($\alpha$) for CVaR-loss/CVaR-regret}\label{alg:alpha-Cal}
\begin{algorithmic}[1]
\State Input: Training dataset $\{\Vzeta_j,\,\Vxi_j\}_{j=1}^{N_{train}}$ and validation dataset $\{\Vzeta_j,\,\Vxi_j\}_{j=1}^{N_{validation}}$ and $\mathcal{A}:=\{\alpha_i\}_{i=1}^{n}\subset[0,\,1]$
\State Train a random forest model $\hat{F}_{\xi|\zeta}$ on $\{\Vzeta_j,\,\Vxi_j\}_{j=1}^{N_{train}}$
\For{$i=1,\dots,n$}
\State //Evaluate $\hat{\x}_i^*(\cdot)$ on empirical distribution of realizations in $\{\Vzeta_j,\,\Vxi_j\}_{j=1}^{N_{validation}}$
\For{$j=1,\dots,N_{validation}$}
\State Solve \eqref{eq:CVaRloss}/\eqref{eq:CVaRregret} with $\alpha_i$ for $\Vzeta:=\Vzeta_j$ in validation set to get optimal $\x_j^*$
\State Let $\hat{\x}_i^*(\Vzeta_j):= \x_j^*$
\EndFor
\State Set $s^i:=\coeffPresc^\alpha_{\hat{F}}\big(\hat{\x}_i^*(\cdot)\big)$ for empirical distribution $\hat{F}$ on $\{\Vzeta_j,\,\Vxi_j\}_{i=1}^{N_{validation}}$
\EndFor
\State Let $i^*:=\argmax_i s^i$ and  set $\alpha^*:=\alpha_{i^*}$ and $\x^*(\cdot):=\hat{\x}_{i^*}^*(\cdot)$
\State Return $(\alpha^*,\x^*(\cdot))$
\end{algorithmic}
\end{algorithm}
}

% \begin{figure*}[ht]
% \vskip 0.2in
% \centerline{\includegraphics[width=\textwidth]{neurIPS.eps}}
% %\centerline{\includegraphics[scale=0.201]{neurIPS.eps}}
% \caption{Out-of-sample performance of (a) relaxed $\x(\cdot)$ (b) binary $\x(\cdot)$}
% \label{fig:OOS}
% \vskip 0.2in
% \end{figure*}

% \begin{figure}[htbp]
% \centerline{\includegraphics[scale=0.5]{Boxplot_relaxed.eps}}
% \caption{Out-of-sample performance of relaxed $\x(\cdot)$}
% \label{fig:OOS}
% \end{figure}

% \begin{figure}[htbp]
% \centerline{\includegraphics[scale=0.5]{Boxplot_binary.eps}}
% \caption{Out-of-sample performance of binary $\x(\cdot)$}
% \label{fig:OOS_bin}
% \end{figure}

Figure \ref{fig:OOS}(a) reports the coefficients of prescriptiveness $\ratioPresc_F(\x^*(\cdot),\hat{\x})$, where $F$ is the test dataset, for the four policies and perturbation levels. More details on the average out-of-sample performance are also presented in Table \ref{table:avgOOS} in Appendix~\ref{appendix:table}. We observe the following: $(i)$ When considering a particular optimization model, the coefficient of prescriptiveness decreases as the magnitude of the distribution shift increases. Indeed, these policies face a more serious robustness challenge as they approach more extreme scenarios beyond what was seen in the train dataset. $(ii)$ When the test set follows the same distribution as the train set, all four policies roughly demonstrate similar performance; however, when this set experiences a distribution shift, DRPCR policies differentiate their performance compared to the alternative ones. $(iii)$ Imposing a more severe distribution shift accentuates this differentiation. For instance, when the mean travel times across the edges are perturbed up to $50\%$ in the test set, DRPCR policies provide a positive coefficient of prescriptiveness, at least over $75\%$ of instances. On the contrary, the alternative policies fail to reach a positive ratio over almost a similar number of instances. This observation is further amplified in the case of $60\%$ perturbation. In this scenario, while CSO, DRCSO, and DRCRO policies fail to return a positive out-of-sample coefficient of prescriptiveness, DRPCR still can reach a positive median of \removed{$2\%$}{\color{black}$4\%$} which can go up to \removed{$11\%$}{\color{black}$18\%$} at its best. 

\pushToAppendix{
\begin{table}[htbp]
    \centering
    \caption{Mean out-of-sample coefficient of prescriptiveness}
    \begin{tabular}{cccccccc}
        \hline
        \multirow{2}{*}{\thead{Problem Type}} & \multirow{2}{*}{\thead{Method}} & \multicolumn{6}{c}{\thead{Level of Perturbation}} \\       
        \cmidrule(lr){3-8}
        & & $0\%$ & $20\%$ & $30\%$  & $40\%$ & $50\%$ & $60\%$\\
        \hline
        
        \multirow{4}{*}{\makecell{Relaxed $\x(\cdot)$}} & \makecell{CSO} & 0.45 & 0.30 & 0.19 & 0.04 & -0.13 & -0.31 \\
        & \makecell{DRCSO} & 0.45 & 0.30 & 0.18 & 0.04 & -0.13 & -0.31 \\
        & \makecell{DRCRO} & 0.45 & 0.30 & 0.18 & 0.04 & -0.13 & -0.32 \\
        & \makecell{DRPCR} & 0.45 & 0.31 & 0.23 & 0.13 & 0.05 & 0.01 \\
        \hline
        
        \multirow{4}{*}{\makecell{Binary $\x(\cdot)$}} & \makecell{CSO} & 0.45 & 0.30 & 0.19 & 0.04 & -0.13 & -0.31 \\
        & \makecell{DRCSO} & 0.44 & 0.30 & 0.19 & 0.06 & -0.09 & -0.25 \\
        & \makecell{DRCRO} & 0.44 & 0.30 & 0.19 & 0.05 & -0.11 & -0.28 \\
        & \makecell{DRPCR} & 0.44 & 0.32 & 0.24 & 0.15 & 0.07 & 0.02\\
        \hline
    \end{tabular}
    \label{tab:mean_oos}
\end{table}
}

\begin{table*}[h]
    \centering
    \caption{Average runtime (minutes) per instance}
    \label{tab:performance}
    \renewcommand{\arraystretch}{1.2}
    \begin{tabular}{@{}l l c c c c c c c@{}}
        \toprule
        \multirow{2}{*}{\thead{Method}} & \multirow{2}{*}{\thead{Type of Problem}} & \multirow{2}{*}{\thead{Algorithm}} & \multicolumn{6}{c}{\thead{Levels of Mean Purturbations ($m$)}}\\
        \cmidrule(lr){4-9}
        & & & 0\% & 20\% & 30\% & 40\% & 50\% & 60\%\\
        \midrule
        DRCSO & relaxed $x(\cdot)$ & \ref{alg:alpha-Cal} & 3.08 & 2.38 & 3.06 & 3.14 & 2.38 & 2.38 \\
        DRCRO & relaxed $x(\cdot)$ & \ref{alg:alpha-Cal} & 2.96 & 2.28 & 2.96 & 3.02 & 2.28 & 2.28 \\
        DRPCR & relaxed $x(\cdot)$ & \ref{alg:alpha-gamma-Cal} & 32.76 & 25.22 & 32.86 & 33.58 & 25.28 & 25.22 \\
        \midrule
        DRCSO & binary $x(\cdot)$ & \ref{alg:alpha-Cal} & 3.08 & 2.64 & 3.04 & 3.42 & 3.02 & 3.40 \\
        DRCRO & binary $x(\cdot)$ & \ref{alg:alpha-Cal} & 3.32 & 3.42 & 3.20 & 3.44 & 3.20 & 3.44 \\
        DRPCR & binary $x(\cdot)$ & \ref{alg:alpha-gamma-Cal} & 34.28 & 38.08 & 33.32 & 37.92 & 33.14 & 37.92 \\
        \bottomrule
    \end{tabular}
\end{table*}

Figure \ref{fig:OOS}(b) depicts the statistics of the 1-norm distance metric of the adaptable policies from the optimal SAA solution. This illustration elucidates the reason behind the superior performance of the DRPCR method compared to others. Indeed, it is notable that DRPCR implicitly utilizes $\bar{x}$ as an anchor for the adaptable policy. In other words, under large distribution shifts, it is able to learn under what context it is worth staying closer to $\bar{x}$, where the relative regret is zero, as explained below Lemma \ref{thm:ratioInterval}.
This phenomenon highlights how the DRPCR applies a completely different form of regularization, compared to prior DRO models.

For further insights, readers are directed to Appendix \ref{app:add_exp}, where an additional set of experiments is presented. In line with the approach in \cite{StochasticOptimizationForests}, the integrality constraint of $\x(\cdot)$ is relaxed in this supplementary investigation. {\color{black} One can refer to Table \ref{tab:performance} for a comparison of runtime of Algorithms \ref{alg:alpha-gamma-Cal} and \ref{alg:alpha-Cal} under both relaxed and binary policies.}

\section{Conclusion}\label{sec:Conclusion}

{\color{black}
The proposed DRPCR model offers an innovative method for calibrating contextual optimization problems and introduces a unique form of regularization, differing from previous DRO models and achieving significantly improved out-of-sample performance. Unfortunately, in its current form, the approach requires considerable training time, i.e. approximately tenfold that of the DRCSO and DRCRO models. This is due to the bisection algorithm needing $\left\lceil{\log_2{(1/\epsilon)}}\right \rceil $ steps to converge. In order to improve tractability (at the expense of optimality), one might consider limiting the admissible policy $\boldsymbol x(\cdot)$ to those with affine dependence on the side information $\boldsymbol \zeta$. The resulting DRPCR takes the form of a smaller optimization problem with size proportional to the number of dimensions of $\boldsymbol \zeta$ rather than $|\Omega_\zeta|$. Affine policies might also facilitate the use of more general non-nested ambiguity sets, which constitute a current limitation of the proposed DRPCR model. Finally, we expect that additional empirical evaluation with other data generating and application environments would certainly benefit our understanding of the value of the presented DRPCR approach.
}

\subsection*{Acknowledgments}
The authors gratefully acknowledge support from the Institut de Valorisation des Données (IVADO), from the Canadian Natural Sciences and Engineering Research Council [RGPIN-2022-05261], and the Canada Research Chair program [CRC-2018-00105].

\bibliography{biblio}
%%%%%%%%%%%%%% APPENDIX %%%%%%%%%%%%%%
\clearpage
\begin{center}
    \LARGE{\textbf{Appendix}}
    \vspace{-0.5em}
    \par\noindent\rule{\textwidth}{1.5pt}
\end{center}
\vspace{0.5cm}

\begin{appendices}
\section{Proofs}\label{appendix:proofs}
\subsection{Proof of Lemma \ref{thm:ratioInterval}}
This follows simply from $\ratioPresc_F(\x(\cdot),
\bar{\x})$ being bound above by 1 for all policy $x(
\cdot)$ and all distribution $F$ due to:
\begin{align*}  
\ratioPresc_F(\x(\cdot),
\bar{\x}) &= 1- \frac{\Expect_F[h(\x(\Vzeta),\V{\xi})]-\Expect_F[\min_{\x'\in\X}h(\x',\V{\xi})]}{\Expect_F[h(\bar{\x},\V{\xi})]-\Expect_F[\min_{\x'\in\X}h(\x',\V{\xi})]} \\
&\leq 1- \frac{\Expect_F[\min_{\x'\in\X}h(\x',\Vxi)]-\Expect_F[\min_{\x'\in\X}h(\x',\V{\xi})]}{\Expect_F[h(\bar{\x},\V{\xi})]-\Expect_F[\min_{\x'\in\X}h(\x',\V{\xi})]} = 1
\end{align*}
when $\Expect_F[h(\xSAA,\Vxi)]-\Expect_F[\min_{\x'\in\X}h(\x',\Vxi)]>0$, and otherwise equal to 1 or $-\infty$ both bounded above by 1. Hence, 
\[\max_{\x(\cdot)\in\policySet} \inf_{F\in\Dist}\ratioPresc_F(\x(\cdot),\bar{\x}) \leq 1.\]
Moreover, if $\bar{\x}\in\policySet$, then we have that
\[\max_{\x(\cdot)\in\policySet} \ratioPresc_\Dist(\x(\cdot),\bar{\x})\geq \ratioPresc_\Dist(\bar{\x},\bar{\x})=\left\{\begin{array}{cl}0&\mbox{if $\Expect_F[h(\bar{\x},\V{\xi})]-\Expect_F[\min_{\x'\in\X}h(\x',\V{\xi})]>0$}\\ 1 & \mbox{ otherwise.}\end{array}\right.\,.\qed\]

\subsection{Proof of Lemma \ref{thm:ERMreduction}}
Let $\tilde{\x}(\cdot)$ be a CSO optimal policy, then necessarily $\tilde{\x}(\cdot)\in\policySet$ since $\tilde{\x}(\Vzeta)\in\X$ for all $\Vzeta$. This confirms that $\tilde{\x}(\cdot)$ is feasible in DRPCR. Next, we can demonstrate optimality through:
    \[\ratioPresc_\Dist(\tilde{\x}(\cdot),\bar{\x})=\ratioPresc_{\bar{F}}(\tilde{\x}(\cdot),\bar{\x})\geq \max_{\x(\cdot)\in\policySet}\ratioPresc_{\bar{F}}(\x(\cdot),\bar{\x}) = \max_{\x(\cdot)\in\policySet}\ratioPresc_\Dist(\x(\cdot),\bar{\x}),\]
since for all $\x(\cdot)\in\policySet$, we have that $\Expect_{F}[h(\x(\Vzeta),\V{\xi})|\Vzeta]\geq \min_{\bm x(\cdot)\in\mathcal{H}} \Expect_{F}[h(\x(\Vzeta),\V{\xi})|\Vzeta]=\Expect_{F}[h(\tilde{\x}(\Vzeta),\V{\xi})|\Vzeta]$ for all $\Vzeta$, which we can show implies  that $\ratioPresc_F(\tilde{\x}(\cdot),\bar{\x})\geq \ratioPresc_F(\x(\cdot),\bar{\x})$. More specifically, if $\Expect_F[h(\xSAA,\Vxi)]=\Expect_F[\min_{\x'\in\X}h(\x',\Vxi)]$, then either $\Expect_{F}[h(\x(\Vzeta),\V{\xi})]=\Expect_F[\min_{\x'\in\X}h(\x',\Vxi)]$ thus 
    \[\Expect_F[\min_{\x'\in\X}h(\x',\Vxi)]=\Expect_{F}[h(\tilde{\x}(\Vzeta),\V{\xi})]\leq  \Expect_{F}[h(\x(\Vzeta),\V{\xi})]=\Expect_F[\min_{\x'\in\X}h(\x',\Vxi)]\]
    meaning that $\ratioPresc_F(\tilde{\x}(\cdot),\bar{\x})=\ratioPresc_F(\x(\cdot),\bar{\x})=1$, or $\Expect_{F}[h(\x(\Vzeta),\V{\xi})]>\Expect_F[\min_{\x'\in\X}h(\x',\Vxi)]$ thus $\ratioPresc_F(\tilde{\x}(\cdot),\bar{\x})\geq -\infty=\ratioPresc_F(\x(\cdot),\bar{\x})$. Alternatively, the case where $\Expect_F[h(\xSAA,\Vxi)]>\Expect_F[\min_{\x'\in\X}h(\x',\Vxi)]$ is straightforward as the function 
    \[f(y):=1- \frac{y-\Expect_F[\min_{\x'\in\X}h(\x',\V{\xi})]}{\Expect_F[h(\bar{\x},\V{\xi})]-\Expect_F[\min_{\x'\in\X}h(\x',\V{\xi})]}\]
    is strictly decreasing.\qed

\subsection{Proof of Proposition \ref{thm:epigraph}}
We first present the DRPCR in epigraph form:
\begin{subequations}
\begin{align}    
    \max_{\gamma,\,\x(\cdot)\in\mathcal{H}} \;\;\;& \gamma\\
    \subto \;\;\;& \ratioPresc_F(\x(\cdot),\bar{\x})\geq \gamma ,\;\forall F\in\Dist\label{eq:prfEpi:C1}\\    
    &0\leq\gamma\leq 1
\end{align}
\end{subequations}
where we added the redundant constraint $\gamma\in[0,\,1]$ since Lemma \ref{thm:ratioInterval} ensures that the optimal value of DRPCR is in this interval.

Focusing on constraint \eqref{eq:prfEpi:C1}, we can then consider two cases for the definition of $\ratioPresc_F(\x(\cdot),\bar{\x})$. In the case that $\Expect_F[h(\bar{\x},\V{\xi})]-\Expect_F[\min_{\x'\in\X}h(\x',\V{\xi})]>0$, one can multiply both sides of the inequality to equivalently obtain:
\[\Expect_{{F}}[h(\x(\Vzeta),\V{\xi})]-\Expect_{{F}}[\min_{\x'\in\X}h(\x',\V{\xi})]\leq (1-\gamma)\left(\Expect_{{F}}[h(\bar{\x},\V{\xi})]-\Expect_{{F}}[\min_{\x'\in\X}h(\x',\V{\xi})]\right)\]
which is equivalent, when rearranging the terms, to:
\begin{equation}
\Expect_{{F}}[h(\x(\Vzeta),\V{\xi})-(1-\gamma) h(\bar{\x},\V{\xi})-\gamma\min_{\x'\in\X}h(\x',\V{\xi})]\leq 0.\label{eq:thmEpi:Cond1}
\end{equation}
In the second case where $\Expect_F[h(\bar{\x},\V{\xi})]=\Expect_F[\min_{\x'\in\X}h(\x',\V{\xi})]$, then constraint \eqref{eq:prfEpi:C1} is equivalent to:
\[\Expect_F[h(\x(\Vzeta),\V{\xi})]=\Expect_F[\min_{\x'\in\X}h(\x',\V{\xi})] \quad\&\quad \gamma\leq 1,\]
yet $\gamma\leq 1$ is redundant while the former condition can equivalently be posed as \eqref{eq:thmEpi:Cond1}.  We are left with
\[\Expect_{{F}}[h(\x(\Vzeta),\V{\xi})-(1-\gamma) h(\bar{\x},\V{\xi})-\gamma\min_{\x'\in\X}h(\x',\V{\xi})]\leq 0\,,\;\forall\,F\in\Dist,\]
which can equivalently be described by $Q(\x(\cdot),\gamma)\leq 0$. One can further conclude that $Q(\x(\cdot),\gamma)\leq 0$ is convex and non-decreasing in $\gamma$ given that it is the supremum of a set of affine non-decreasing functions:
\[Q(\x(\cdot),\gamma)=\sup_{F\in\Dist}\Expect_{{F}}[h(\x(\Vzeta),\V{\xi})-h(\bar{\x},\V{\xi})]+\gamma(\Expect_F[h(\bar{\x},\V{\xi})-\min_{\x'\in\X}h(\x',\V{\xi})]),\]
with $h(\bar{\x},\V{\xi})\geq\min_{\x'\in\X}h(\x',\V{\xi})$ for all $\Vxi$ since $\bar{\x}\in\X$.\qed

\subsection{Proof of Proposition \ref{thm:epigraph_reformulation}}
Letting $g(\x,\Vxi,\gamma):=h(\x,\V{\xi})-\left((1-\gamma) h(\bar{\x},\V{\xi})+\gamma \min_{x'\in\X}h(\x',\V{\xi})\right)$, we have that
\begin{align*}  
\psi(\gamma):=&\min_{\x(\cdot)\in\mathcal{H}}Q(\x(\cdot),\gamma)\\
=&\min_{\x(\cdot)\in\mathcal{H}}\sup_{F\in\DistCVaR(\bar{F},\alpha)}\;\Expect_{{F}}\Big[g(\x(\Vzeta),\V{\xi},\gamma)\Big]\\
=&\min_{\x(\cdot)\in\mathcal{H}}\Expect_{\bar{F}}\Bigg[\CVaR^{\alpha}_{\bar{F}}\Big(g(\x(\Vzeta),\Vxi,\gamma)|\Vzeta\Big)\Bigg]\\
=&\min_{\x(\cdot)\in\mathcal{H}}\Expect_{\bar{F}}\Bigg[\inf_t\;t+\frac{1}{1-\alpha}\Expect_{\bar{F}}\bigg[\max\Big(0,g(\x(\Vzeta),\Vxi,\gamma)-t\Big)|\Vzeta\bigg]\Bigg]\\
=&\Expect_{\bar{F}}\Bigg[\inf_{\x\in\X,t}\;t+\frac{1}{1-\alpha}\Expect_{\bar{F}}\bigg[\max\Big(0,g(\x,\Vxi,\gamma)-t\Big)|\Vzeta\bigg]\Bigg],
\end{align*}
where we exploit the infimum representation of CVaR and the interchangeability property of expected value operators (see \cite{SHAPIRO2017377} and reference therein).
Given that $\bar{F}$ is a discrete distribution as described in Assumption \ref{ass:CVaR}, one can compute $\psi(\gamma)$ by solving for each scenario $\Vzeta_\omega$ with $\omega\in\Omega_\zeta$ the problem \eqref{eq:probScenBased:C2}.
% \begin{subequations}\label{eq:probScenBased}
% \begin{eqnarray}
%     \phi_\omega(\gamma):=\quad\quad\min_{\x\in\X,t,\V{s}}&&t+\frac{1}{1-\alpha}\sum_{\omega'\in\Omega_\xi} \Prob_{\bar{F}}(\Vxi=\Vxi_{\omega'}|\Vzeta=\Vzeta_\omega)s_{\omega'}\\
%     \subto && s_{\omega'} \geq h(\x,\V{\xi}_{\omega'})-\left((1-\gamma) h(\bar{\x},\V{\xi}_{\omega'})+\gamma \min_{x'\in\X}h(\x',\V{\xi}_{\omega'})\right) - t\,,\;\forall \omega'\in\Omega_\xi\label{eq:probScenBased:C2}\\
%     &&s_{\omega'} \geq 0\,,\;\forall \omega'\in\Omega_\xi.
% \end{eqnarray}
% \end{subequations}
Based on the solution of problem \eqref{eq:probScenBased} for each $\omega\in\Omega_\zeta$, one can obtain $\psi(\gamma):=\sum_{\omega\in\Omega_\zeta}\Prob_{\bar{F}}(\Vzeta=\Vzeta_\omega) \phi_\omega(\gamma)$ together with a potentially feasible policy $\x(\Vzeta):=\x_{\omega(\Vzeta)}^*$, 
where $\omega(\Vzeta)=\argmin_{\omega\in\Omega_\zeta}\|\Vzeta-\Vzeta_\omega\|$ and $\x_\omega$ refers to the minimizer of problem \eqref{eq:probScenBased}.

The function $\phi_\omega(\gamma)$ is non-decreasing in $\gamma$ since $\gamma$ only appears in constraint \eqref{eq:probScenBased:C2}, which can be rewritten as:
\[s_{\omega'} \geq h(\x,\V{\xi}_{\omega'})-h(\bar{\x},\V{\xi}_{\omega'})+\left( h(\bar{\x},\V{\xi}_{\omega'}) - \min_{x'\in\X}h(\x',\V{\xi}_{\omega'})\right)\gamma - t\,,\;\forall \omega'\in\Omega_\xi\,.\]
Since the right-hand side of this constraint is non-decreasing in $\gamma$, due to $h(\bar{\x},\V{\xi}_{\omega'}) \geq \min_{x'\in\X}h(\x',\V{\xi}_{\omega'})$, one can concludes that the minimum of \eqref{eq:probScenBased} cannot decrease when $\gamma$ is increased, since the feasible set is reduced.

We further note that problem \eqref{eq:probScenBased} can be reduced to a linear program when $\X$ is polyhedral and $h(\x,\Vxi_{\omega'})$ is linear programming representable for all $ \omega'\in\Omega_\xi$. 
For example, in the context of a portfolio optimization, where $\X$ is the probability simplex and $h(\x,\Vxi):=-\Vxi^T\x$, we have that problem \eqref{eq:thmEpi:mainProb} reduces to:
\begin{eqnarray*}
    \max_{\gamma,\{\x^\omega,t^\omega,\V{s}^\omega\}_{\omega\in\Omega_\zeta}}&&\gamma\\
    \text{subject to}&&\sum_{\omega\in\Omega_\zeta} \Prob_{\bar{F}}(\Vzeta=\Vzeta_\omega) \left(t^\omega+\frac{1}{1-\alpha}\sum_{\omega'\in\Omega_\xi} \Prob_{\bar{F}}(\V{\xi}=\V{\xi}_{\omega'}|\Vzeta=\Vzeta_\omega)s_{\omega'}^\omega\right)\leq 0\\
    &&s_{\omega'}^\omega \geq \V{\xi}_{\omega'}^T\x^\omega-(1-\gamma) \V{\xi}_{\omega'}^T\bar{\x}-\gamma \min_{x'\in\X}\V{\xi}_{\omega'}^T\x' - t^\omega\,,\;\forall \omega'\in\Omega_\xi,\,\omega\in\Omega_\zeta\\
    &&s^\omega \geq 0 \,,\;\forall \omega\in\Omega_\zeta\\
    && \sum_{i=1}^{n_x} \x_i^\omega = 1,\;\forall \omega\in\Omega_\zeta\\
    && \x^\omega \geq 0,\;\forall \omega\in\Omega_\zeta\\
    && 0\leq \gamma\leq 1\,.
\end{eqnarray*}
Alternatively, in the context of a shortest path problem (see Section \ref{sec:Experiments} for details), we have that problem \eqref{eq:thmEpi:mainProb} reduces to a mixed integer linear program:
\begin{eqnarray*}
    \max_{\gamma,\{\x^\omega,t^\omega,\V{s}^\omega\}_{\omega\in\Omega_\zeta}}&&\gamma\\
    \text{subject to}&&\sum_{\omega\in\Omega_\zeta} \Prob_{\bar{F}}(\Vzeta=\Vzeta_\omega) \left(t^\omega+\frac{1}{1-\alpha}\sum_{\omega'\in\Omega_\xi} \Prob_{\bar{F}}(\V{\xi}=\V{\xi}_{\omega'}|\Vzeta=\Vzeta_\omega)s_{\omega'}^\omega\right)\leq 0\\
    &&s_{\omega'}^\omega \geq \V{\xi}_{\omega'}^T\x^\omega-(1-\gamma) \V{\xi}_{\omega'}^T\bar{\x}-\gamma \min_{x'\in\X}\V{\xi}_{\omega'}^T\x' - t^\omega\,,\;\forall \omega'\in\Omega_\xi,\,\omega\in\Omega_\zeta\\
    &&s^\omega \geq 0 \,,\;\forall \omega\in\Omega_\zeta\\
&&\sum_{j:(o,j)\in\mathcal{A}} x_{(o,j)}^\omega-\sum_{j:(j,o)\in\mathcal{A}} x_{(j,o)}^\omega=1,\;\forall \omega\in\Omega_\zeta\\
&&\sum_{j:(d,j)\in\mathcal{A}} x_{(d,j)}^\omega-\sum_{j:(j,d)\in\mathcal{A}} x_{(j,d)}^\omega=-1,\;\forall \omega\in\Omega_\zeta\\
&&\sum_{j:(i,j)\in\mathcal{A}} x_{(i,j)}^\omega-\sum_{j:(j,i)\in\mathcal{A}} x_{(j,i)}^\omega=0,\;\forall i\in\mathcal{V}\,\backslash\,\{o,\,d\},\,\omega\in\Omega_\zeta\\
&&\x_{(i,j)}^\omega\in\{0,\,1\},\;\forall (i,j)\in\mathcal{A},\,\omega\in\Omega_\zeta\\
    && 0\leq \gamma\leq 1\,.\hfill\qed
\end{eqnarray*}

{\color{black}
\subsection{Proof of Lemma \ref{convergence_algorithm1}}
\EDmodified{%From Proposition~\ref{thm:epigraph} we know that $\phi_\omega(\gamma)$ and $\psi(\gamma):=\sum_{\omega\in\Omega_\zeta} \Prob_{\bar{F}}(\Vzeta=\Vzeta_\omega) \phi_\omega(\gamma)$ are non-decreasing. 
One can first easily verify that in Algorithm \ref{alg:bisection}, we have that $\Delta:=\gamma^+-\gamma^-$ is initially equal to 1 and reduces by a factor of 2 at every iteration. The algorithm therefore necessarily terminates after $\lceil \log_2(1/\epsilon)  \rceil$ iterations. When $\X$ is polyhedral and $h(\x,\V{\xi})$ is linear programming representable, problem \eqref{nesterCVaRproblem} reduces to a linear program that can be solved in polynomial time with respect to $|\Omega_\xi|$, $n_x$, $n_\xi$, the size of the LP representation of $\mathcal{X}$, and of $h(\boldsymbol x,\boldsymbol\xi)$ (see \cite{ellipsoidMethod} and \cite{polyTimeLP}). Given that this problem is solved $|\Omega_\zeta|$ at each iteration of the algorithm. We conclude that the total run time of the algorithm is polynomial with respect to all of these quantities. \qed
}
%From Proposition~\ref{thm:epigraph} we know that $Q(\x(\cdot),\gamma)$ is a convex non-decreasing function of $\gamma$.  Thus we know from \cite{boyd2004convex} (Chapter 4.2.5) that the bisection algorithm will converge in exactly $\lceil \log_2(1/\epsilon)  \rceil$ iterations. Furthermore, it is easy to see that problem \eqref{eq:probScenBased} is polynomially solvable with its complexity  being affected by the number of scenarios   $|\Omega_\zeta|$ and $|\Omega_\xi|$ as well as the complexity of the $h(\bm x,\bm \xi)$ and $\mathcal{X}$. Hence, the algorithm will terminates in polynomial time with respect to $\log(1/\epsilon).$
}

\subsection{\EDmodified{Proof of Proposition~\ref{generalized_proposition2}} \label{sec:appgenCVaR}}

The proof follows similar steps as Proposition~\ref{thm:epigraph_reformulation}. Let $g(\x,\Vxi,\gamma):=h(\x,\V{\xi})-\left((1-\gamma) h(\bar{\x},\V{\xi})+\gamma \min_{x'\in\X}h(\x',\V{\xi})\right)$, we have that
\begin{align}  
\psi(\gamma)=&\min_{\x(\cdot)\in\mathcal{H}}\sup_{F\in\DistCVaR(\bar{F},\alpha)}\;\Expect_{{F}}\Big[g(\x(\Vzeta),\V{\xi},\gamma)\Big]\notag\\
=&\min_{\x(\cdot)\in\mathcal{H}}\sup_{\V{p}:\V{p}\geq 0,\sum_{\omega\in\Omega_\zeta}p_\omega=1, d_\zeta(\V{p},\bar{\V{p}})\leq r_\zeta}\sum_{\omega\in\Omega_\zeta} p_\omega \sup_{\V{q}:\V{q}\geq 0,\sum_{\omega'\in\Omega_\xi}q_{\omega'}=1, d_\xi(\V{q},\bar{\V{q}}^\omega)\leq r_\xi} \sum_{\omega'\in\Omega_\xi} q_{\omega'} g(\x(\Vzeta_\omega),\Vxi_{\omega'},\gamma)\notag\\
=& \sup_{\V{p}:\V{p}\geq 0,\sum_{\omega\in\Omega_\zeta}p_\omega=1, d_\zeta(\V{p},\bar{\V{p}})\leq r_\zeta}\sum_{\omega\in\Omega_\zeta} p_\omega \min_{\x\in\X} \sup_{\V{q}:\V{q}\geq 0,\sum_{\omega'\in\Omega_\xi}q_{\omega'}=1, d_\xi(\V{q},\bar{\V{q}}^\omega)\leq r_\xi} \sum_{\omega'\in\Omega_\xi} q_{\omega'}g(\x(\Vzeta_\omega),\Vxi_{\omega'},\gamma)\notag\\
=& \sup_{\V{p}:\V{p}\geq 0,\sum_{\omega\in\Omega_\zeta}p_\omega=1, d_\zeta(\V{p},\bar{\V{p}})\leq r_\zeta}\sum_{\omega\in\Omega_\zeta} p_\omega \bar{\phi}_\omega(\gamma),\label{eq:genConst10b}
\end{align}
with
\begin{align}
\bar{\phi}_\omega(\gamma):=
\min_{\x\in\X} \sup_{\V{q}:\V{q}\geq 0,\sum_{\omega'\in\Omega_\xi}q_{\omega'}=1, d_\xi(\V{q},\bar{\V{q}}^\omega)\leq r_\xi} \sum_{\omega'\in\Omega_\xi}  q_{\omega'}g(\x(\Vzeta_\omega),\Vxi_{\omega'},\gamma)    .\label{eq:phibardef}
\end{align}

Denoting $\V{g}\in\Re^{|\Omega'|}$ with $g_{\omega'}:=g(\x(\Vzeta_\omega),\Vxi_{\omega'},\gamma)$, one can derive a reformulation of the inner supremum in \eqref{eq:phibardef} as an infimum following:
\begin{align*}
\sup_{\V{q}:\V{q}\geq 0,\sum_{\omega'\in\Omega_\xi}q_{\omega'}=1, d_\xi(\V{q},\bar{\V{q}}^\omega)\leq r_\xi}  \V{q}^T\V{g}&=\sup_{\V{q}} \inf_{\V{\lambda}\geq 0,t,\alpha\geq 0}  \V{q}^T \V{g} +\V{\lambda}^T\V{q} +t(1-\1^T\V{q}) + \alpha(r_\xi-d_\xi(\V{q},\bar{\V{q}}^\omega))\\
&=\inf_{\V{\lambda}\geq 0,t,\alpha\geq 0} \sup_{\V{q}} \V{q}^T \V{g} +\V{\lambda}^T\V{q} +t(1-\1^T\V{q}) + \alpha(r_\xi-d_\xi(\V{q},\bar{\V{q}}^\omega))\\
&=\inf_{\V{\lambda}\geq 0,t,\alpha\geq 0} t +r_\xi\alpha + d_*(\V{g}+\V{\lambda}-t,\alpha,\bar{\V{q}}^\omega),
\end{align*}
where $d_*(\V{v},\alpha,\bar{\V{q}}^\omega):=\sup_{\V{q}} \V{v}^T\V{q}-\alpha d_\xi(\V{q},\bar{\V{q}}^\omega)$ is the perspective of the convex conjugate  of $d_\xi(\V{q},\bar{\V{q}}^\omega)$.

By joining this reformulation with the minimization in $\x$, we get %that problem \eqref{nesterCVaRproblem} should be replaced with:
\begin{subequations}\label{eq:nestergeneralCVaRproblem}
\begin{eqnarray}
    \min_{\x\in\X,t,\alpha\geq 0,\V{s}}&& t+r_\xi\alpha + d_*(\V{s},\alpha,\bar{\V{q}}^\omega)\\
    \text{subject to}&&s_{\omega'}\geq g(\x(\Vzeta_\omega),\Vxi_{\omega'},\gamma)-t,\;\forall \omega'\in\Omega'
\end{eqnarray}    
\end{subequations}
which concludes the proof.\qed

\section{Acceleration strategy for Algorithm \ref{alg:bisection}}\label{Appendix:Bisection}

One can possibly accelerate the convergence rate on the bisection Algorithm \ref{alg:bisection} by exploiting the fact that $\psi(\cdot)$ is a convex function when $\X$ is convex. Indeed, for the current interval $[\gamma^-,\gamma^+]$, $\psi(\gamma)$ can be under- and over-estimated, see Figure~\ref{fig:bisection} (right). The procedure can be described as follows. First, we construct a line that will underestimate $\psi$   by identifying a subgradient of the function at $\tilde\gamma$. This can be computed analytically since 
\begin{equation*}
\begin{array}{rl}
\psi(\gamma)&:=\displaystyle\Expect_{\bar F}\left[\min_{x\in\mathcal{X}} \CVaR_\alpha\left(h(\x,\V{\xi})-\Big((1-\gamma) h(\bar{\x},\V{\xi})+\gamma \min_{x'\in\X}h(\x',\V{\xi})\Big)\middle|\zeta\right)\right]\\[2ex]

&=\displaystyle\Expect_{\bar F}\left[\min_{x\in\mathcal{X}} \sup_{F\in\DistCVaR(\bar{F},\alpha)}\Expect_F\left[h(\x,\V{\xi})-\Big((1-\gamma) h(\bar{\x},\V{\xi})+\gamma \min_{x'\in\X}h(\x',\V{\xi})\Big)\middle|\zeta\right]\right]\\[2ex]

&\geq\displaystyle\Expect_{\bar F}\left[ \sup_{F\in\DistCVaR(\bar{F},\alpha)}\min_{x\in\mathcal{X}}\Expect_F\left[h(\x,\V{\xi})-\Big((1-\gamma) h(\bar{\x},\V{\xi})+\gamma \min_{x'\in\X}h(\x',\V{\xi})\Big)\middle|\zeta\right]\right]\\   
  
&\geq\displaystyle\Expect_{\bar F}\left[ \min_{x\in\mathcal{X}}\Expect_{ F^*_{\xi|\zeta}}\left[h(\x,\V{\xi})-\Big((1-\gamma) h(\bar{\x},\V{\xi})+\gamma \min_{x'\in\X}h(\x',\V{\xi})\Big)\middle|\Vzeta\right]\right]\\   

&=\displaystyle\underbrace{\Expect_{\bar F}\left[ \min_{x\in\mathcal{X}}\Expect_{ F^*_{\xi|\zeta}}\left[h(\x,\V{\xi})-  h(\bar{\x},\V{\xi}) \right]\right]}_{\mymbox{offset \quoteIt{$a$}}} + 
\gamma\underbrace{\Expect_{\bar F}\left[\Expect_{ F^*_{\xi|\zeta}}\left[ h(\bar{\x},\V{\xi}) - \min_{x'\in\X}h(\x',\V{\xi}) \right]\right]}_{\mymbox{slope \quoteIt{$b$}}},   
\end{array}
\end{equation*}
where $F_{\xi|\zeta}^*$ is the conditional probability given $\Vzeta$ of any member (hopefully a maximizer) of $\DistCVaR(\bar{F},\alpha)$. Note that the first inequality is tight  based on 
Sion's minimax theorem (see \cite{sion58:minimax}) given that $\DistCVaR(\bar{F},\alpha)$ is compact, while the second is tight as long as $F_{\xi|\zeta}^*$ achieves the supremum. Such a maximizer can be identified using:
% \[F^*_{\xi|\zeta}\in\argmax_{F_{\xi|\zeta}\in\mathcal{M}(\Omega_\zeta):\Prob_{F_{\xi|\zeta}}(\V{\xi})\leq (1-\alpha)^{-1}\Prob_{\bar{F}}(\Vxi|\Vzeta),\,\forall\Vxi}\Expect_{F_{\xi|\zeta}
% }\left[h(\x^*_\gamma(\bm\zeta),\V{\xi})-\Big((1-\gamma) h(\bar{\x},\V{\xi})+\gamma \min_{x'\in\X}h(\x',\V{\xi})\Big)\right]\]
\[F^*_{\xi|\zeta}\in\argmax_{\tiny{
\begin{array}{cc}
F_{\xi|\zeta}\in\mathcal{M}(\Omega_\zeta):\\
\Prob_{F_{\xi|\zeta}}(\V{\xi})\leq(1-\alpha)^{-1}\Prob_{\bar{F}}(\V{\xi}|\Vzeta),\,\forall\V{\xi}
\end{array}}}\Expect_{F_{\xi|\zeta}
}\left[h(\x^*_\gamma(\bm\zeta),\V{\xi})-\Big((1-\gamma) h(\bar{\x},\V{\xi})+\gamma \min_{x'\in\X}h(\x',\V{\xi})\Big)\right]\]
where $\x^*_\gamma(\bm\zeta)$ is the minimizer of \eqref{eq:probScenBased} with $\Vzeta_\omega=\Vzeta$ since $(\bm x^*_\gamma(\cdot),F^*)$, with $F^*$ as the composition of $\bar{F}$ marginalized on $\Vzeta$ and $F^*_{\Vxi|\Vzeta}$,\footnote{Namely, $P_{F^*}(\Vxi)=P_{\bar{F}}(\Vxi)$ and $P_{F^*}(\Vxi|\Vzeta)=P_{F^*_{\xi|\zeta}}(\Vxi)$ for all $\Vzeta$.} is a saddle point of:
\[g(\x(\cdot),F):=\Expect_{{F}}\Big[h(\x(\Vzeta),\V{\xi})-\Big((1-\gamma) h(\bar{\x},\V{\xi})+\gamma \min_{x'\in\X}h(\x',\V{\xi})\Big)\Big].\]
Such a $F^*_{\Vxi|\Vzeta}$ can be obtained as a side product of solving problem \eqref{eq:probScenBased} using the optimal dual variables associated with constraint \eqref{eq:probScenBased:C2}. If we denote by $\gamma_u := a/b$ then the right bound of the interval can be updated to $\gamma^{+'} := \min(\gamma^+,\gamma_u)$. 

\begin{figure}[ht!]
\centerline{\includegraphics[scale=0.55]{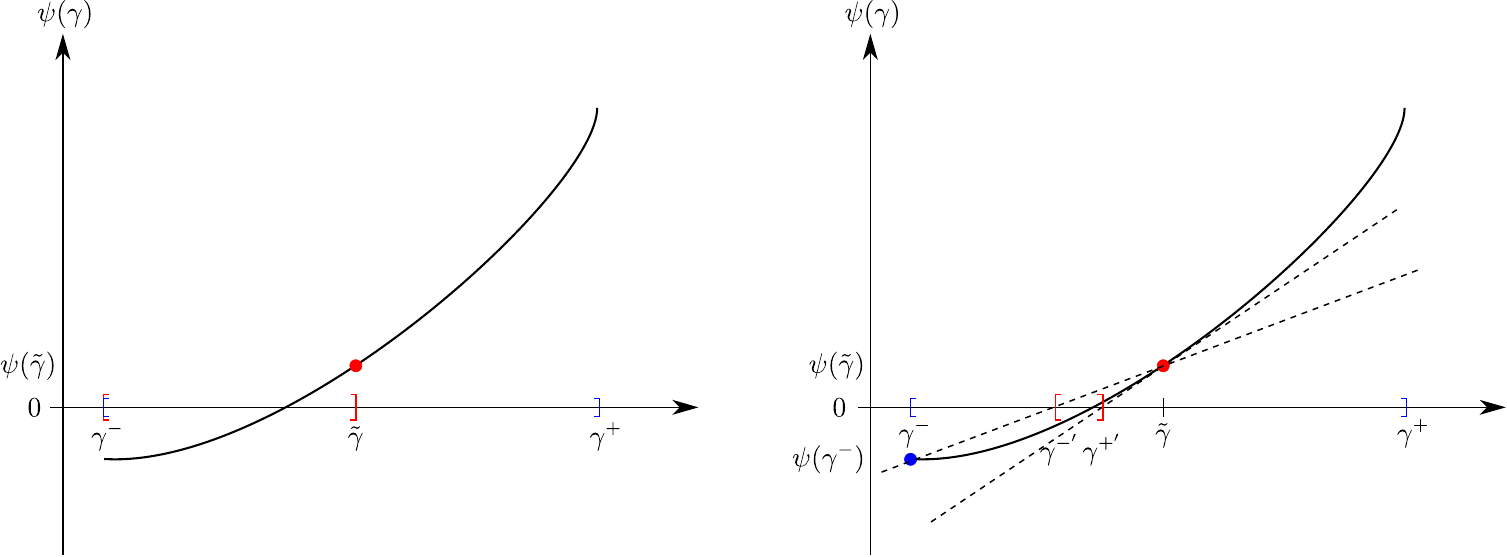}}
\caption{Visualization of the basic (left) and accelerated (right) bisection algorithm. The blue squared brackets indicate the current estimated interval containing the optimal $\gamma^*$ and the red squared brackets indicate the interval in the next iterations. The right graph also visualizes the over and under estimators of $\psi(\gamma)$.
}
\label{fig:bisection}
\end{figure}

The second step is to construct an overestimator. If $\psi(\tilde\gamma) > 0$, then we evaluate $\psi(\gamma^{-})$ and construct the line that passes through $(\gamma^{-},\psi(\gamma^{-}))$ and $(\tilde\gamma,\psi(\tilde\gamma))$. If $\psi(\tilde\gamma) < 0$ then  we evaluate $\psi(\gamma^{+})$ and construct the line that passes through $(\gamma^{+},\psi(\gamma^{+}))$ and $(\tilde\gamma,\psi(\tilde\gamma))$. We denote the point for which the line evaluates to zero as $\gamma_o$, and update the left bound of the interval to $\gamma^{-'} := \max(\gamma^-,\gamma_o)$. Hence, the new interval is given by $[\gamma^{-'},\gamma^{+'}]\subseteq[\gamma^{-},\gamma^{+}]$, which would potentially significantly reduce the search space.

We conclude this section by commenting that the accelerated bisection algorithm could require up to two evaluations of the $\psi$ function at each iteration instead of a single one as described in the original algorithm.

\newpage
\section{Algorithms for calibrating the size of the ambiguity sets}\label{Appendix:Algo}

\begin{algorithm}[htbp]
\caption{Algorithm for calibrating the size of the ambiguity set ($\alpha$) for DRPCR}\label{alg:alpha-gamma-Cal}
\begin{algorithmic}[1]
\State Input: Training dataset $\{\Vzeta_j,\,\Vxi_j\}_{j=1}^{N_{train}}$ and validation dataset $\{\Vzeta_j,\,\Vxi_j\}_{j=1}^{N_{validation}}$ and $\mathcal{A}:=\{\alpha_i\}_{i=1}^{n}\subset[0,\,1]$
\State Train a random forest model $\hat{F}_{\xi|\zeta}$ on $\{\Vzeta_j,\,\Vxi_j\}_{j=1}^{N_{train}}$
\State Let $\hatFC$ be the composition of $\hat{F}_{\xi|\zeta}$ with empirical distribution $\hat{F}_\zeta$ of $\Vzeta$ in the training set $\{\Vzeta_j\}_{j=1}^{N_{train}}$
\For{$i=1,\dots,n$}
\State //Construct $\hat{\x}_i^*(\cdot)$ with $\alpha_i$ and $\hatFC$
\State Solve DRPCR with $\alpha_i$ and $\hatFC$ to get $\gamma_i^*$ \label{step:DRPCRsolve}
\State //Evaluate $\hat{\x}_i^*(\cdot)$ on empirical distribution of realizations in $\{\Vzeta_j,\,\Vxi_j\}_{j=1}^{N_{validation}}$
\For{$j=1,\dots,N_{validation}$}
\State Solve \eqref{eq:probScenBased} with $\gamma_i^*$, $\alpha_i$, and replacing $\Prob_{\bar{F}}(\Vxi=\Vxi_{\omega'}|\Vzeta=\Vzeta_\omega)$ with $\Prob_{\hat{F}_{\xi|\zeta_j}}(\Vxi_{\omega'})$ to get optimal $\x_j^*$
\State Let $\hat{\x}_i^*(\Vzeta_j):= \x_j^*$
\EndFor
\State Set $s^i:=\coeffPresc^\alpha_{\hat{F}}\big(\hat{\x}_i^*(\cdot)\big)$ for empirical distribution $\hat{F}$ on $\{\Vzeta_j,\,\Vxi_j\}_{i=1}^{N_{validation}}$
\EndFor
\State Let $i^*:=\argmax_i s^i$ and  set $\alpha^*:=\alpha_{i^*}$, $\gamma^*:=\gamma_{i^*}$, and $\x^*(\cdot):=\hat{\x}_{i^*}^*(\cdot)$
\State Return $(\alpha^*,\gamma^*,\x^*(\cdot))$
\end{algorithmic}
\end{algorithm}

\begin{algorithm}[htbp!]
\caption{Algorithm for calibrating the size of the ambiguity set ($\alpha$) for CVaR-loss/CVaR-regret}\label{alg:alpha-Cal}
\begin{algorithmic}[1]
\State Input: Training dataset $\{\Vzeta_j,\,\Vxi_j\}_{j=1}^{N_{train}}$ and validation dataset $\{\Vzeta_j,\,\Vxi_j\}_{j=1}^{N_{validation}}$ and $\mathcal{A}:=\{\alpha_i\}_{i=1}^{n}\subset[0,\,1]$
\State Train a random forest model $\hat{F}_{\xi|\zeta}$ on $\{\Vzeta_j,\,\Vxi_j\}_{j=1}^{N_{train}}$
\For{$i=1,\dots,n$}
\State //Evaluate $\hat{\x}_i^*(\cdot)$ on empirical distribution of realizations in $\{\Vzeta_j,\,\Vxi_j\}_{j=1}^{N_{validation}}$
\For{$j=1,\dots,N_{validation}$}
\State Solve \eqref{eq:CVaRloss}/\eqref{eq:CVaRregret} with $\alpha_i$ for $\Vzeta:=\Vzeta_j$ in validation set to get optimal $\x_j^*$
\State Let $\hat{\x}_i^*(\Vzeta_j):= \x_j^*$
\EndFor
\State Set $s^i:=\coeffPresc^\alpha_{\hat{F}}\big(\hat{\x}_i^*(\cdot)\big)$ for empirical distribution $\hat{F}$ on $\{\Vzeta_j,\,\Vxi_j\}_{i=1}^{N_{validation}}$
\EndFor
\State Let $i^*:=\argmax_i s^i$ and  set $\alpha^*:=\alpha_{i^*}$ and $\x^*(\cdot):=\hat{\x}_{i^*}^*(\cdot)$
\State Return $(\alpha^*,\x^*(\cdot))$
\end{algorithmic}
\end{algorithm}

\section{Generation of random covariance matrix with arbitrary variances}\label{CovMat}
A random covariance matrix for the random vector of $(\Vzeta,\V{\xi})$ is generated based on a two-step procedure that follows. The first step consists in generating a random symmetric positive-definite matrix described in Algorithm \ref{alg:SPD}, a method implemented in the sklearn.datasets.make\_spd\_matrix function of scikit-learn machine learning library in Python.   

\begin{algorithm}[H]
\caption{Algorithm for generating random symmetric positive-definite matrix}\label{alg:SPD}
\begin{algorithmic}[1]
\State Input: Dimension of the square matrix $n_\Vzeta + n_\V{\xi}$
\State Generate random square matrix $A_{n_\Vzeta + n_\V{\xi}}$ sampling from the uniform distribution $\mathcal{U}_{[0,1]}$
\State Construct the symmetric matrix $M=A^\top A$
\State Decompose $M$ with Singular Value Decomposition (SVD) method as $M = U \Sigma V^\top$ 
\State Generate random diagonal matrix $S$ sampling from the uniform distribution $\mathcal{U}_{[0,1]}$
\State Construct $\Sigma^{'} = S + J$ where $J$ is the square matrix of ones with the size of $n_\Vzeta + n_\V{\xi}$
\State Get the symmetric positive-definite matrix as $M' = U \Sigma^{'} V^\top$ 
\State Return $M'$
\end{algorithmic}
\end{algorithm}

Given the vector of standard deviations for $(\Vzeta,\V{\xi})$ denoted by $[\V{\sigma}^\top_\Vzeta\,\,\V{\sigma}^\top_\V{\xi}]^\top$ and also a random symmetric positive-definite matrix generated by Algorithm \ref{alg:SPD}, one can implement the second stage described in Algorithm \ref{alg:Covariance} to get a random covariance matrix with arbitrary standard deviations of $[\V{\sigma}^\top_\Vzeta\,\,\V{\sigma}^\top_\V{\xi}]^\top$.

\begin{algorithm}[H]
\caption{Algorithm for generating random covariance matrix with arbitrary standard deviations}\label{alg:Covariance}
\begin{algorithmic}[1]
\State Input: Random symmetric positive-definite matrix ($M$) and vector of standard deviations $[\V{\sigma}^\top_\Vzeta\,\,\V{\sigma}^\top_\V{\xi}]^\top$
\State Convert matrix $M$ into its associated correlation matrix $\mbox{Corr}=\Big(\mbox{diag(M)}\Big)^{-\frac{1}{2}}M\Big(\mbox{diag(M)}\Big)^{-\frac{1}{2}}$
\State Get the arbitrary covariance matrix of $\mbox{Cov}=\mbox{diag}\Big(\Big[\begin{array}{c} \V{\sigma}_\Vzeta \\ \V{\sigma}_\V{\xi} \end{array}\Big]\Big)\,\,\big(\mbox{Corr}\big)\,\,\mbox{diag}\Big(\Big[\begin{array}{c} \V{\sigma}_\Vzeta \\ \V{\sigma}_\V{\xi} \end{array}\Big]\Big)$
\State Return Cov
\end{algorithmic}
\end{algorithm}

\section{Average out-of-sample coefficient of prescriptiveness}\label{appendix:table}
\begin{table}[htbp]
    \centering
    \caption{Average out-of-sample coefficient of prescriptiveness}\label{table:avgOOS}
    \begin{tabular}{cccccccc}
        \hline
        \multirow{2}{*}{\thead{Problem Type}} & \multirow{2}{*}{\thead{Method}} & \multicolumn{6}{c}{\thead{Level of Perturbation}} \\       
        \cmidrule(lr){3-8}
        & & $0\%$ & $20\%$ & $30\%$  & $40\%$ & $50\%$ & $60\%$\\
        \hline
        
        \multirow{4}{*}{\makecell{Relaxed $\x(\cdot)$}} & \makecell{CSO} & 0.45 & 0.30 & 0.19 & 0.04 & -0.13 & -0.31 \\
        & \makecell{DRCSO} & 0.45 & 0.30 & 0.18 & 0.04 & -0.13 & -0.31 \\
        & \makecell{DRCRO} & 0.45 & 0.30 & 0.18 & 0.04 & -0.13 & -0.32 \\
        & \makecell{DRPCR} & 0.45 & 0.31 & 0.23 & 0.13 & 0.05 & 0.01 \\
        \hline
        
        \multirow{4}{*}{\makecell{Binary $\x(\cdot)$}} & \makecell{CSO} & 0.45 & 0.30 & 0.19 & 0.04 & -0.13 & -0.31 \\
        & \makecell{DRCSO} & 0.44 & 0.30 & 0.19 & 0.06 & -0.09 & -0.25 \\
        & \makecell{DRCRO} & 0.44 & 0.30 & 0.19 & 0.05 & -0.11 & -0.28 \\
        & \makecell{DRPCR} & 0.44 & 0.32 & 0.24 & 0.15 & 0.07 & 0.02\\
        \hline
    \end{tabular}
    \label{tab:mean_oos}
\end{table}

\section{{\color{black}Additional experiments}}\label{app:add_exp}

\begin{figure}[htbp!]
\centerline{\includegraphics[width=\textwidth]{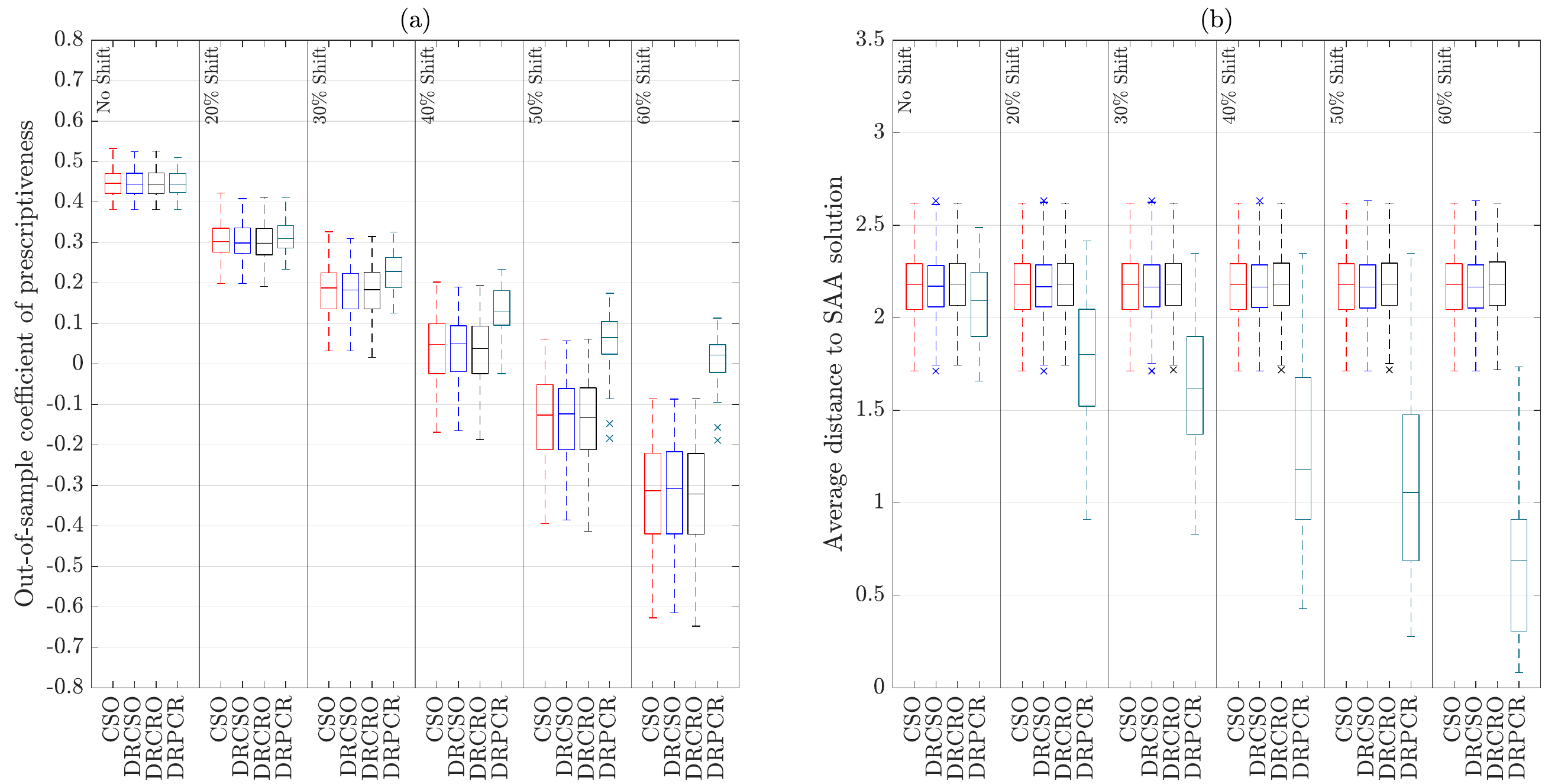}}
\caption{Shortest path problem (relaxed version): (a) statistics of the out-of-sample coefficient of prescriptiveness {\color{black} (lower values indicate worse performance).} (b) statistics of $\mathbb{E}_{\EDmodified{\breve{F}}}[\|\boldsymbol{x^*}(\boldsymbol \zeta) - \boldsymbol{\hat x}\|_1]$ where \EDmodified{$\breve F$ is the} out-of-sample distribution {\color{black} (lower values reflect a closer proximity to the SAA solution).}}
\label{fig:OOS_relaxed}
\end{figure}

While the experiments in Section \ref{sec:Experiments} consider an exact version of the shortest path problem, to be closer to the setting proposed in \cite{StochasticOptimizationForests}, we also conduct a second set of experiments where $\x(\cdot)$ represents relaxed variables. Figure \ref{fig:OOS_relaxed} (a) illustrates the coefficients of prescriptiveness obtained from the optimal relaxed policies. These results, in general, are aligned with the ones spotted in Figure \ref{fig:OOS} (a); however, one remarks the following. Firstly, the results derived from CSO remain exactly the same as the binary case. This stems from the fact that optimal relaxed CSO decisions are known to be integral for the stochastic shortest path problems; conversely, this is not the case for DRCSO, DRCRO, and DRPCR where robustness breaks the linearity of the objective. Secondly, comparing Figures \ref{fig:OOS} (a) and \ref{fig:OOS_relaxed} (a) reveals that forcing DRCSO and DRCRO to propose binary policies enhances their out-of-sample performance, surpassing those of CSO. Indeed, this setting seems to provide these two approaches the chance to better prepare for potential distribution shifts; however, despite their enhanced performance, the highest degree of robustness to distribution shift remains associated with DRPCR policies. Thirdly, this comparative analysis yields counter-intuitive empirical evidence that out-of-sample performance might be slightly improved when imposing integrality constraints on the three robust models. We hypothesize that this might be caused by the additional flexibility of the relaxed models, which makes them more susceptible to overfitting their assumed stochastic models.

\end{appendices}
\end{document}